\documentclass[EJP]{ejpecp} 

\SHORTTITLE{GSDEs with Reflecting Boundary Conditions} 

\TITLE{Stochastic Differential Equations\\ Driven by $G$-Brownian Motion\\ with Reflecting Boundary Conditions} 

\AUTHORS{%
  Yiqing~Lin\footnote{Shandong University, China. \EMAIL{yiqing.lin@univ-rennes1.fr}} \footnote{Universit\'{e} de Rennes 1, France.}
    }


\AMSSUBJ{60H10} 

\SUBMITTED{June 29, 2011} 
\ACCEPTED{December 29, 2012} 

\VOLUME{18}
\YEAR{2013}
\PAPERNUM{0}
\DOI{vVOL-PID}

\ABSTRACT{In this paper, we introduce the idea of stochastic integrals with respect to an increasing process in the $G$-framework and extend $G$-It\^o's formula. Moreover, we study the solvability of the scalar valued stochastic differential equations driven by $G$-Brownian motion with reflecting boundary conditions (RGSDEs).}

\newcounter{newlist} 
\newenvironment{mylist}[1][]{
\begin{list}{\textnormal{(\arabic {newlist})}} 
    { 
    \usecounter{newlist} 
     \setlength{\labelsep}{0.5em}
     \setlength{\leftmargin}{0.5em} 
     \setlength{\rightmargin}{0cm} 
     \setlength{\parsep}{0pt}
     \setlength{\itemsep}{3pt} 
     \setlength{\itemindent}{0em} 
     \setlength{\listparindent}{2em} 
    }}
{\end{list}}
\newcounter{nnnewlist} 
\newenvironment{mmmylist}[1][]{
\begin{list}{\textnormal{(\arabic {nnnewlist})}} 
    { 
    \usecounter{nnnewlist} 
     \setlength{\labelsep}{0.5em} 
     \setlength{\leftmargin}{1.7em} 
     \setlength{\rightmargin}{0cm} 
     \setlength{\parsep}{0pt} 
     \setlength{\itemsep}{3pt} 
     \setlength{\itemindent}{0em} 
     \setlength{\listparindent}{2em} 
    }}
{\end{list}}
\newcounter{nelist} 
\newenvironment{mlist}[1][]{
\begin{list}{\textnormal{(A\arabic {nelist})}} 
    { 
    \usecounter{nelist} 
     \setlength{\labelsep}{0.5em} 
     \setlength{\leftmargin}{2.5em} 
     \setlength{\rightmargin}{0cm} 
     \setlength{\parsep}{0pt} 
     \setlength{\itemsep}{3pt} 
     \setlength{\itemindent}{0em} 
     \setlength{\listparindent}{2em} 
    }}
{\end{list}}

\begin{document}



\section{Introduction}
\noindent In the classical framework, Skorokhod \cite{S1,S2} first introduced diffusion processes with reflecting boundaries in the 1960s. Since then, reflected solutions to stochastic differential equations (SDEs) and Backward SDEs (BSDEs) have been investigated by many authors. For the one-dimensional case, El Karoui \cite{EL}, El Karoui and Chaleyat-Maurel \cite{EC} and Yamada \cite{Y} studied reflected SDEs (RSDEs) on a half-line and El Karoui et al. \cite{EKPPQ} obtained the solvability of reflected BSDEs. For the multidimensional case, the existence of weak solutions to reflected SDEs on a smooth domain was proved by Stroock and Varadhan \cite{SV1}. Subsequently, Tanaka \cite{T} solved the similar problem on a convex domain by a direct approach based on the solution to the Skorokhod problem. Furthermore, Lions and Sznitman \cite{LS} extended these results to a non-convex domain. The corresponding results for reflected BSDEs can be found in Gegout-Petit and Pardoux \cite{GP}, Ramasubramanian \cite{R} and Hu and Tang \cite{HT} and others. \\[6pt]
Motivated by uncertainty problems, risk measures and super-hedging in finance, Peng \cite {P2, P3} introduced a framework of time consistent nonlinear expectation $\mathbb{E}[\cdot]$, i.e., $G$-expectation, in which a new type of Brownian motion was constructed and the corresponding stochastic calculus was established. In order to solve the super-replication problem in an uncertainty volatility model, Denis and Martini \cite{DM} independently introduced a notion of upper expectation and the related capacity theory. Moreover, a stochastic integral of It\^o's type under a class of non-dominated probability measures was formulated. Recently, Hu et al. \cite{DHP} found there is a strong link that connects these two frameworks, that is, the $G$-expectation $\mathbb{E}[\cdot]$ can be represented by a concrete weakly compact family $\mathcal{P}_G$ of probability measures:
\[
\mathbb{E}[X]=\sup_{\mathbb{P}\in\mathcal{P}_G}E^\mathbb{P}[X],\ X\in L^1_G(\Omega).
\]
Then, a Choquet capacity $\bar{C}(\cdot)$ can be naturally introduced to the $G$-framework:
\[
\bar{C}(A):=\sup_{\mathbb{P}\in\mathcal{P}_G}\mathbb{P}(A),\ A\in\mathcal{B}(\Omega),
\]
by which we can have the following definition to the concept of ``quasi-surely'', similar to the one in Denis and Martini \cite{DM}: A set $A\subset\Omega$ is polar if $\bar{C}(A)=0$; and a property holds ``quasi-surely'' (q.s. for short) if it holds outside a polar set. In these two frameworks, a stochastic integral of It\^o type is defined following a usual procedure, that is, giving a definition first for some simple integrands and then completing the spaces of integrands under the norm induced by the upper expectation related to $\mathcal{P}_G$. This norm is much stronger than that in the classical case and thus, the space of integrands is smaller than the classical one. In other words, some additional regularity assumption should be imposed on the integrands to ensure that the integrals are well defined. Using these notions of stochastic calculus in the $G$-framework, the existence and uniqueness results for some types of SDEs driven by $G$-Brownian motion (GSDEs) can be obtained (cf. Peng \cite {P3}, Gao \cite {G} and Lin and Bai \cite{BL}). For the reason stated above, the authors who studied GSDEs always assumed the following condition on the coefficients of the equations: for each $x\in\mathbb{R}$,
\[
f_\cdot(x),\ g_\cdot(x)\in M^2_G([0, T]).
\]
At this price, all results in the works for GSDEs listed above hold in the ``quasi-surely'' (q.s.) sense, i.e., outside a polar set, and all the processes are immediately aggregated.\\[6pt]
\noindent Closely related to the $G$-framework, Soner et al. \cite{STZ1, STZ3, STZ2} have established another type of ``quasi-sure'' stochastic analysis and also a complete theory for second order BSDEs (2BSDEs) under a uniform Lipschitz condition on the coefficients. In that framework, another notion of ``quasi-surely'' was issued, which means that a property holds $\mathbb{P}$-a.s., for each probability measure $\mathbb{P}\in\mathcal{P}_H$, which is a class of local martingale measure.
Obviously, this definition of ``quasi-surely'' is weaker than the one made by $G$-capacity. In this weaker sense, we can consider the stochastic integral with respect to the canonical $B$ under each probability measure $\mathbb{P}\in\mathcal{P}_H$, respectively and we only need that these integrands meet the requirement for formulating a stochastic integral with respect to a local martingale. Thus, this type of setting for 2BSDEs ensures that we can treat the case that the coefficients have less regularity but that all the properties can only hold $\mathbb{P}$-a.s., for each $\mathbb{P}\in\mathcal{P}_H$.
Following the pioneering work of Soner et al. \cite{STZ2}, Matoussi et al. \cite{MPZ1} have studied the problem of reflected 2BSDEs with a lower obstacle.\\[6pt]
\noindent The aim of this paper is to study
the solvability of stochastic differential
equations driven by $G$-Brownian motion with reflecting boundary conditions (RGSDEs) in the sense of ``quasi-surely'' defined by Denis et al. \cite{DHP}. The scalar valued RGSDE that we consider is defined as following:
\begin{equation}\label{ABC}
\left\{\begin{aligned}
&X_t=x+\int^t_0 f_s(X_s)ds+\int^t_0 h_s(X_s)d\langle B \rangle_s +\int^t_0 g_s(X_s) dB_s +K_t, \ 0\leq t\leq T,\ q.s.; \\
&X_t\geq S_t,\ 0\leq t\leq T,\ q.s.;\ \int^T_0 (X_t-S_t)dK_t=0,\ q.s.,
\end{aligned}\right.
\end{equation}
where $\langle B\rangle$ is the quadratic variation process of $G$-Brownian motion $B$ and $K$ is an increasing process that pushes the solution $X$ upwards to remain above the obstacle $S$ in a minimal way. Similarly to how the uniqueness results for classical reflected SDEs have been proved, the corresponding ones for RGSDEs can also be deduced from a priori estimates. Moreover, a solution in $M^p_G([0, T])$ to (\ref{ABC}) can be constructed 
by fixed-point iteration. Because of the reason that we have already explained, we need in addition to some assumption on the coefficients $f$, $h$ and $g$, which is similar to that in Peng \cite {P3}, Gao \cite {G} and Lin and Bai \cite{BL}, a regularity assumption on $S$ to ensure that $K$ stays in the space $M^p_G([0, T])$. To establish the comparison theorem, we need to develop an extension of $G$-It\^{o}'s formula to deal with such a process $X$, which involves both the stochastic integrals and an increasing process. This extended $G$-It\^o's formula can have its own interest and may be used in other situations.\\[6pt]
This paper is organized as follows: Section 2 introduces notation and results in the $G$-framework which are necessary for what follows. Section 3 introduces the stochastic calculus with respect to an increasing process in the $G$-framework. Section 4 studies reflected $G$-Brownian motion and Section 5 presents our main results.
\section{$G$-Brownian motion, $G$-capacity and $G$-stochastic calculus}
\noindent The main purpose of this section is to recall some preliminary results in the $G$-framework, which are necessary later in the text. The reader interested
in a more detailed description of these notions is referred to
Denis et al. \cite{DHP}, Gao \cite{G} and Peng \cite{P3}.
\subsection{$G$-Brownian motion}\noindent Adapting the approach in Peng \cite{P3}, let $\Omega$ be the space of all $\mathbb{R}$-valued
continuous paths with $\omega_0 =0$ 
equipped
with the distance
\[
\rho(\omega^1, \omega^2):=\sum^\infty_{N=1} 2^{-N} ((\max_
{0\leq t\leq N} | \omega^1_t-\omega^2_t|)  \wedge 1),
\]
$B$ the canonical process and $C_{l,Lip}
(\mathbb{R}^{n})$ the collection of all local Lipschitz functions on $\mathbb{R}^{n}$. For a fixed $T\geq0$, the space of finite dimensional cylinder random variables is defined by  
\[
L_{ip}^0 (\Omega_T):=\{\varphi(B_{t_1}, \ldots, B_{t_n}):n\geq
1,\ 0\leq t_1\leq \ldots\leq t_n \leq T,\ \varphi \in C_{l,Lip}
(\mathbb{R}^{n})\},
\]
on which $\mathbb{E}[\cdot]$
is a sublinear functional that satisfies: for all $X$, $Y\in L^0_{ip}(\Omega_T)$,
\begin{mmmylist}
\item\textbf{Monotonicity:} if $X\geq Y$, then $\mathbb{E}[X]\geq\mathbb{E}[Y];$
\item\textbf{Sub-additivity}:\ $\mathbb{E}[X]-\mathbb{E}[Y]\leq\mathbb{E}[X-Y]$;
\item\textbf{Positive homogeneity:}\ $\mathbb{E}[\lambda
X]=\lambda\mathbb{E}[X]$, for all $\lambda\geq 0$;
\item\textbf{Constant translatability:}\ $\mathbb{E}[X+c]=\mathbb{E}[X]+c$, for all $c\in\mathbb{R}$.
\end{mmmylist}
The triple $(\Omega, L_{ip}^0 (\Omega_T), \mathbb{E})$ is called a sublinear expectation space.
\begin{definition}\label{de 2.1.4}
A scalar valued random variable $X\in L_{ip}^0 (\Omega_T)$ is $G$-normal distributed with parameters $(0, [\underline{\sigma}^2, \overline{\sigma}^2])$, i.e., $X\sim\mathcal{N}(0, [\underline{\sigma}^2, \overline{\sigma}^2])$, 
if for each $\varphi \in C_{l,Lip} (\mathbb{R})$, $u(t, x):=\mathbb{E}[\varphi(x+\sqrt{t}X)]$
is a viscosity solution to the following PDE on $\mathbb{R}^+\times\mathbb{R}$:
\[
\left\{\begin{aligned}
&\frac{\partial u}{\partial t}-G\bigg(\frac{\partial^2 u}{\partial x^2}\bigg)=0;\\
&u|_{t=0} = \varphi,
\end{aligned}\right.
\]
where
\[
G(a):=\frac{1}{2} (a^+\overline{\sigma}^2-a^-\underline{\sigma}^2),\ a\in\mathbb{R}.
\]
\end{definition}
\begin{remark}
Without loss of generality, we always assume that $\overline{\sigma}^2=1$ in what follows.
\end{remark}
\begin{definition}\label{de 2.1.5}
We call a sublinear expectation $\mathbb{E}: L_{ip}^0(\Omega_T) \rightarrow \mathbb{R}$ a $G$-expectation if the canonical process $B$ is a $G$-Brownian motion under
$\mathbb{E}[\cdot]$, that is, for each $0\leq s\leq t\leq T$, the increment $B_t - B_s \sim \mathcal{N}(0, [(t-s)\underline{\sigma}^2, (t-s)])$ and for all $n>0$, $0\leq
t_1\leq\ldots\leq t_n\leq T$ and $\varphi\in
C_{l,Lip}(\mathbb{R}^{n})$,
\[
\mathbb{E}[\varphi(B_{t_1},\ldots,B_{t_{n-1}},B_{t_n}-B_{t_{n-1}})]=\mathbb{E}[\psi(B_{t_1},\ldots,B_{t_{n-1}})],
\]
where
$\psi(x_1,\ldots,x_{n-1}):=\mathbb{E}[\varphi(x_1,\ldots,x_{n-1},\sqrt{t_n-t_{n-1}}B_1)]$.
\end{definition}
\noindent For $p\geq 1$, we denote by $L^p_G(\Omega_T)$ the completion of $L_{ip}^0(\Omega_T)$ under the Banach norm $\mathbb{E}[|\cdot|^p]^\frac{1}{p}$.
\subsection{$G$-capacity} \noindent Derived in Denis et al. \cite{DHP}, $G$-expectation $\mathbb{E}[\cdot]$ can  be viewed as an upper expectation $\bar{\mathbb{E}}[\cdot]$ associated with a weakly compact family $\mathcal{P}_G$ of probability measures on $L^1_G(\Omega_T)$, i.e.,
\[\mathbb{E}[X]=\bar{\mathbb{E}}[X]:=\sup_{\mathbb{P}\in\mathcal{P}_G}E^\mathbb{P}[X],\ X\in L^1_G(\Omega_T).\]
In this sense,  the domain of $G$-expectation can be extended from $L^1_G(\Omega_T)$ to the space of all $\mathcal{B}(\Omega_T)$ measurable random variables $L^0(\Omega_T)$ by setting
\[
\bar{\mathbb{E}}[X]:=\sup_{\mathbb{P}\in\mathcal{P}_G}E^\mathbb{P}[X],\ X\in L^0(\Omega_T).
\]
\noindent Naturally, we can define a corresponding regular Choquet capacity on $\Omega$:
\[
\bar{C}(A):=\sup_{\mathbb{P}\in\mathcal{P}_G}\mathbb{P}(A),\ A\in\mathcal{B}(\Omega),
\]
with respect to which, we have the following notions:
\begin{definition}\label{de 2.1.6}
A set $A\in\mathcal{B}(\Omega)$ is called polar if $\bar{C}(A)=0.$ A property
is said to hold quasi-surely (q.s.) if it holds outside a polar
set.
\end{definition}
\begin{definition}\label{qcon}
A random variable $X$ is said to be quasi-continuous (q.c.) if for any arbitrarily small $\varepsilon>0$, there exists an open set $O\subset \Omega$ with $\bar{C}(O)<\varepsilon$ such that $X$ is continuous in $\omega$ on $O^c$.
\end{definition}
\begin{definition}
We say that a random variable $X$ has a q.c. version if there exists a q.c. random variable $Y$ such that $X=Y$, q.s..
\end{definition}
\noindent In the language of $G$-capacity, Denis et al. \cite{DHP} proved that for each $p\geq 1$, the function space $L^p_G(\Omega_T)$ has a dual representation, which is much more explicit to verify:
\begin{theorem}\label{lpg}
\[
L^p_G(\Omega_T)=\{X\in L^0(\Omega_T): X\ has\ a\ q.c.\ version,\ \lim_{N\rightarrow+\infty}\bar{\mathbb{E}}[|X|^p{\mathbb{1}}_{|X|>N}]=0\}.
\]
\end{theorem}
\noindent Unlike in the classical framework, the downwards monotone convergence theorem only holds true for a sequence of random variables from a subset of $L^0(\Omega_T)$ (cf. Theorem 31 in Denis et al. \cite{DHP}).
\begin{theorem}\label{CTD}
Let $\{X^n\}_{n\in\mathbb{N}}\subset L_G^1(\Omega_T)$ be such that $X^n\downarrow X$, q.s., then $\bar{\mathbb{E}}[X^n]\downarrow\bar{\mathbb{E}}[X]$.
\end{theorem}
\begin{remark}
We note that dominated convergence theorem does not exist in the $G$-framework, even though we assume that $\{X^n\}_{n\in\mathbb{N}}$ is a sequence in $L^1_G(\Omega_T)$. The lack of this theorem is one of the main difficulties we shall overcome in the following sections.
\end{remark}
\subsection{$G$-stochastic calculus} \noindent In Peng \cite{P3}, generalized It$\hat{o}$ integrals with respect to $G$-Brownian motion are established:
\begin{definition}\label{caocao}
A partition of $[0,T]$ is a finite
ordered subset $\pi^N_{[0, T]}=\{t_0, t_1,\ldots,t_N\}$ such that
$0=t_0<t_1<\ldots<t_N=T$. We set
\[\mu(\pi^N_{[0,T]}):=\max_{k=0, 1,\ldots,N-1}|t_{k+1}-t_k|.\]
For each $p\geq1$, define
\[
M_G^{p,0}([0,T]):=\bigg\{ \eta_t=\sum_{k=0}^{N-1}\xi_{k}{\mathbb{1}}_{[t_{k},t_{k+1})}(t):
\xi_{k}\in L_G^p(\Omega_{t_{k}})\bigg\},
\]
and denote by $M_G^p([0,T])$ the completion of
$M_G^{p,0}([0,T])$ under the norm
\[||\eta||_{M_G^p([0,T])}:=\bigg(\frac{1}{T}\int_0^T\bar{\mathbb{E}}[|\eta_t|^p]dt\bigg)^\frac{1}{p}.\]
\end{definition}
\begin{remark}\label{remconna}
By Definition \ref{caocao}, 
if $\eta$ is an element in $M^p_G([0,T])$, then there exists a sequence of processes $\{\eta^n\}_{n\in\mathbb{N}}$ in $M_G^{p,0}([0,T])$, such that $
\lim\limits_{n\rightarrow+\infty}\int_0^T\bar{\mathbb{E}}[|\eta^n_t-\eta_t|^p]dt\rightarrow 0
$. It is easily observed that for almost every $t\in[0,T]$, $\{\eta^n_t\}_{n\in\mathbb{N}}\subset L^p_G(\Omega_t)$ and $\bar{\mathbb{E}}[|\eta^n_t-\eta_t|^p]\rightarrow 0$, thus $\eta_t$ is an element in $L^p_G(\Omega_t)$.
\end{remark}
\begin{definition} For each $\eta\in M_G^{2,0}([0,T])$, we define
\[
\mathcal{I}_{[0, T]}(\eta)=\int_0^T\eta_sdB_s
:=\sum_{k=0}^{N-1}\xi_k(B_{t_{k+1}}-B_{t_k}).
\]
The mapping $\mathcal{I}_{[0, T]}: M_G^{2, 0}([0,T])\rightarrow L^2_G(\Omega_T)$ is continuous and linear and thus, can be uniquely extended to $\mathcal{I}_{[0, T]}: M_G^2([0,T])\rightarrow L^2_G(\Omega_T)$. Then, for each $\eta\in M_G^2([0,T])$, the stochastic integral with respect to $G$-Brownian motion $B$ is defined by $\int_0^T\eta_sdB_s:=\mathcal{I}_{[0, T]}(\eta)$.
\end{definition}
\noindent Unlike the classical theory, the quadratic variation process of $G$-Brownian motion $B$ is not always a deterministic process (unless $\overline{\sigma}=\underline{\sigma}$) and it can be formulated in $L^2_G(\Omega_t)$ by 
\[\langle B\rangle_t:=\lim_{\mu(\pi^N_{[0,t]})\rightarrow 0}\sum^{N-1}_{k=0}(B_{t^n_{k+1}}-B_{t^n_k})^2=B^2_t-2\int^t_0B_sdB_s.\]
\begin{definition} For each $\eta\in M_G^{1,0}([0,T])$, we define
\[
\mathcal{Q}_{[0,T]}(\eta)=\int_0^T\eta_sd\langle B\rangle_s
:=\sum_{k=0}^{N-1}\xi_k(\langle B\rangle_{t_{k+1}}-\langle B\rangle_{t_k}).
\]
The mapping $\mathcal{Q}_{[0,T]}: M_G^{1, 0}([0,T])\rightarrow L^1_G(\Omega_T)$ is continuous and linear and thus, can be uniquely extended to $\mathcal{Q}_{[0, T]}: M_G^1([0,T])\rightarrow L^1_G(\Omega_T)$. Then, for each $\eta\in M_G^1([0,T])$, the stochastic integral with respect to the quadratic variation process $\langle B\rangle$ is defined by $\int_0^T\eta_sd\langle B\rangle_s:=\mathcal{Q}_{[0,T]}(\eta)$.
\end{definition}
\noindent In view of the dual formulation of $G$-expectation, as well as the properties of the quadratic variation process $\langle B\rangle$ in the $G$-framework, the following BDG type inequalities are obvious.
\begin{lemma}\label{th 2.1.2}
Let $p\geq 1$, $\eta\in M^p_G([0,T])$ and $0\leq s\leq t\leq T$. Then, 
\[
\bar{\mathbb{E}}\bigg[\sup_{s\leq u\leq t}\bigg|\int^u_s\eta_rd\langle B\rangle_r\bigg|^p\bigg]\leq |t-s|^{p-1}\int_t^s\bar{\mathbb{E}}[|\eta_u|^p]du.
\]
\end{lemma}
\begin{lemma}\label{th 2.1.1}
Let $p\geq 2$, $\eta\in M^p_G([0,T])$ and $0\leq s\leq t\leq T$. Then, 
\[
\bar{\mathbb{E}}\bigg[\sup_{s\leq u\leq t}\bigg|\int^u_s\eta_rdB_r\bigg|^p\bigg]\leq C_p\bar{\mathbb{E}}\bigg[\bigg|\int_s^t|\eta_u|^2du\bigg|^\frac{p}{2}\bigg]\leq C_p|t-s|^{\frac{p}{2}-1}\int_s^t\bar{\mathbb{E}}[|\eta_u|^p]du,
\]
where $C_p$ is a positive constant independent of $\eta$.
\end{lemma}
\section{Stochastic calculus with respect to an increasing process}
\noindent In this section, we define the stochastic integrals with respect to an increasing
process with continuous paths, and then we extend $G$-It\^o's
formula to the case where an increasing process appears in the dynamics. In the following, $C$ and $M$ denote two positive constants whose values may vary from line to line. 
\subsection{Stochastic integrals with respect to an increasing process}
\begin{definition}
We denote by $M_c([0, T])$ the collection of all q.s. continuous processes $X$ whose paths $X_\cdot(\omega): t\mapsto X_t(\omega)$ are continuous in $t$ on $[0, T]$ outside a polar set $A$. 
\end{definition}
\begin{remark}\label{diffcon}For example, from the proofs of Theorem 2.1 and Theorem 2.2 in Gao \cite{G}, $(\int^t_0\eta_sdB_s)_{0\leq t\leq T}$ and $(\int^t_0\eta_sd\langle B\rangle_s)_{0\leq t\leq T}$ have continuous modifications in $M_c([0, T])$.
\end{remark}
\begin{definition}
We denote by $M_I([0, T])$ the collection of q.s. increasing processes $K\in M_c([0, T])$
whose paths $K_\cdot(\omega): t\mapsto K_t(\omega)$ are increasing in $t$ on $[0, T]$ outside a polar set $A$.
\end{definition}
\begin{remark}Obviously, an increasing process $K$ in $M_I([0, T])$ has q.s. finite total variation
on $[0, T]$ and thus, its quadratic variation is q.s. 0.
\end{remark}
\begin{definition}\label{increa}
We define, for a fixed $X\in M_c([0, T])$, the stochastic integral with respect to a given $K\in M_I([0, T])$ by
 \begin{equation}\label{II}
\bigg(\int^T_0 X_t dK_t\bigg)(\omega)
=\left\{
\begin{aligned}
\int^T_0 X_t(\omega)dK_t(\omega)\ &,\
\omega\in A^c;\\
0\ &,\ \omega\in A,
\end{aligned}\right.
\end{equation}
where $A$ is a polar set and on the complementary of which, $X_\cdot(\omega)$ is continuous and $K_\cdot(\omega)$ is continuous and increasing in $t$.
\end{definition} 
\begin{remark}
Because for a fixed $\omega\in A^c$, the function $X_\cdot(\omega)$ is continuous and the function $K_\cdot(\omega)$ is of bounded variation on $[0, T]$, the Riemann-Stieltjes integral on the right-hand side always exists (cf. Hildebrandt \cite{HI}). Thus, (\ref{II}) is well defined. Similar definitions can be made for those $X$ whose paths are q.s. piecewisely continuous and without discontinuity of the second kind , i.e., for each $\omega\in A^c$, the function $X_\cdot(\omega)$ is discontinuous at a finite number of points and these discontinuous points are removable or of the first kind.
\end{remark}
\begin{remark}\label{sequen}
Given a sequence of refining partitions $\{\pi^N_{[0,T]}\}_{N\in\mathbb{N}}$ (i.e., $\pi^N_{[0,T]}\subset \pi^{N+1}_{[0,T]}$, for all $N\in\mathbb{N}$) such that $\mu(\pi^N_{[0,T]})\rightarrow 0$, we set a sequence of binary functions:
\begin{equation}\label{constr}
\mathcal{V}_{[0,T]}^N(X, K)(\omega)
:=\sum^{N-1}_{k=0}
X_{u^N_k}(\omega)(K_{t^N_{k+1}}(\omega)-K_{t^N_k}(\omega)),
\end{equation}
where $u^N_k\in[t^N_k, t^N_{k+1})$. For a fixed $\omega\in A^c$, by the Heine-Cantor theorem, $X_\cdot(\omega)$ and $K_\cdot(\omega)$ are uniformly continuous in $t$ on $[0, T]$. Therefore, we can find an $M_\omega>0$ such that $K_T(\omega)<M_\omega$, then, for any arbitrarily small $\varepsilon>0$, there exists a $\delta>0$ such that for all
$|t-s|<\delta$, $|X_t(\omega)-X_s(\omega|)<\varepsilon/M_\omega$. It is sufficient to choose an $N_0\in\mathbb{N}$ such that $\mu(\pi^{N_0}_{[0,T]})<\delta$, then, for all $N>N_0$, 
\[
\bigg|\mathcal{V}_{[0,T]}^N(X, K)(\omega)-\bigg(\int^T_0X_tdK_t\bigg)(\omega)\bigg|<\varepsilon,
\]
from which we deduce
\begin{equation}\label{icon}
\mathcal{V}_{[0,T]}^N(X, K)\rightarrow \int^T_0X_tdK_t,\ q.s.,\ as\ N\rightarrow +\infty.
\end{equation}
The construction of sequence (\ref{constr}) provides a q.s. approximation to the stochastic integral $\int^T_0X_tdK_t$. We note that the convergence (\ref{icon}) depends only on the sequence of refined partitions $(\pi^N_{[0,T]})_{N\in\mathbb{N}}$ but is independent of the selection of the points of division and the representatives $X_{u^N_k}$ on $[t^N_k, t^N_{k+1})$,
$k=0, 1,\ldots, N-1$, $N\in\mathbb{N}$.
\end{remark}
\noindent The following propositions can be verified directly by Definition \ref{increa} and the Heine-Cantor theorem.
\begin{proposition}\label{propo1}
Let $X$, $X^1$, $X^2\in M_c([0, T])$, $K$, $K^1$, $K^2\in M_I([0, T])$ and $0\leq s\leq r \leq t\leq T$, then we have
\begin{mylist}
\item $\int^t_sX_udK_u=\int^r_sX_udK_u+\int^t_rX_udK_u$, q.s.;
\item $\int^t_s(\alpha X^1_u+X^2_u)dK_u=\alpha\int^t_sX^1_udK_u+\int^t_sX^2_udK_u$, q.s., where $\alpha\in L^0(\Omega_s)$;
\item\label{sum} $\int^t_sX_ud(K^1\pm K^2)_u=\int^t_sX_udK^1_u \pm \int^t_sX_udK^2_u$, q.s..
\end{mylist}
\end{proposition}
\begin{remark}
By a classical argument, a q.s. continuous and bounded variation process can be viewed as the difference of two increasing processes $K_1-K_2$, where $K_1$ and $K_2\in M_I([0, T])$. By Proposition \ref{propo1} (\ref{sum}), the stochastic integral with respect to $K_1-K_2$ can be defined in the same way as Definition \ref{increa}. 
\end{remark}
\begin{proposition}
Let $X\in M_c([0, T])$ and $K\in M_I([0, T])$, 
then the integral $\int^\cdot_0X_sdK_s$ is q.s. continuous in $t$, i.e., $(\int^t_0X_sdK_s)_{0\leq t\leq T}\in M_c([0, T])$.
\end{proposition}
\noindent As shown above, (\ref{II}) defines a random variable $\int^T_0X_tdK_t$ in $L^0(\Omega_T)$. A natural question arises: if we assume that for some appropriate $p$ and $q$, $X\in M^p_G([0, T])$ and $K\in M^q_G([0, T])$, can this random variable $\int^T_0X_tdK_t$ be verified as an element in $L^1_G(\Omega_T)$ or not? In general, the answer is negative. This is because the integrability of $X$ and $K$ cannot ensure the quasi-continuity of $\int^T_0X_tdK_t$ (cf. Definition \ref{qcon} and Theorem \ref{lpg}). More precisely, the pathwise convergence (\ref{icon}) is not necessarily uniform in $\omega$ outside a polar set $A$ and it is hard to verify directly the convergence in the sense of $L^1_G(\Omega_T)$ due to the lack of the dominated convergence theorem in the $G$-framework. However, in some special cases, a proper sequence $\{\mathcal{V}_{[0,T]}^N(X, K)\}_{N\in\mathbb{N}}$ approximating to $\int^T_0X_tdK_t$ can be found and thus, the quasi-continuity is inherited during the approximation.
\begin{proposition}
Let $K\in M_I([0, T])\cap M^2_G([0, T])$, $K_T\in L^2_G(\Omega_T)$ and $\phi: \mathbb{R}\rightarrow \mathbb{R}$ is a Lipschitz function, then $\int^T_0\phi(K_t)dK_t$ is an element in $L^1_G(\Omega_T)$. 
\end{proposition}
\noindent\textbf{Proof:} Consider a sequence of refining partitions $\{\pi^N_{[0,T]}\}_{N\in\mathbb{N}}$ mentioned in Remark \ref{sequen} and define the sequence of approximation: for each $N\in\mathbb{N}$, 
\[
\mathcal{V}_{[0,T]}^N(\phi(K), K)(\omega)
=\sum^{N-1}_{k=0}
\phi(K_{t^N_k})(\omega)(K_{t^N_{k+1}}(\omega)-K_{t^N_k}(\omega)).
\]
From the explanation in Remark \ref{remconna}, 
we can always assume that at the points of division, $K_{t^N_k}\in L^2_G(\Omega_T)$, $k =0,  1,\ldots, N-1$, $N\in \mathbb{N}$.
As $K$ is increasing, we have
\begin{align*}
\bigg|\mathcal{V}_{[0,T]}^N(\phi(K), K)-\int^T_0 \phi(K_t)dK_t\bigg|
&\leq\bigg|\int^T_0 \bigg(\sum^{N-1}_{k=0}|K_{t^N_{k+1}}-K_{t^N_k}|{\mathbb{1}}_{[t^N_k, t^N_{k+1})}(t)\bigg) dK_t\bigg|\\
&\leq\sum^{N-1}_{k=0}|K_{t^N_{k+1}}-K_{t^N_k}|^2 \downarrow 0,\ q.s.,\ \textnormal{as}\ N\rightarrow +\infty.
\end{align*}
On the other hand, it is easy to verify by Theorem \ref{lpg} that for all $N\in\mathbb{N}$, $\mathcal{V}_{[0,T]}^N(\phi(K), K)$ and $\sum^{N-1}_{k=0}|K_{t^N_{k+1}}-K_{t^N_k}|^2\in L^1_G(\Omega_T)$. Then, by Theorem \ref{CTD}, we obtain
\[
\bar{\mathbb{E}}\bigg[\bigg|\mathcal{V}_{[0,T]}^N(\phi(K), K)-\int^T_0 \phi(K_t)dK_t\bigg|\bigg]
\leq \bar{\mathbb{E}}\bigg[\sum^{N-1}_{k=0}|K_{t^N_{k+1}}-K_{t^N_k}|^2\bigg]\downarrow0,\ \textnormal{as}\ N\rightarrow +\infty.
\]
From the completeness of $L^1_G(\Omega_T)$ under $\bar{\mathbb{E}}[|\cdot|]$, 
we deduce the desired result.\hfill{}$\square$
\begin{remark}
To verify that for all $N\in \mathbb{N}$, $\mathcal{V}_{[0,T]}^N(\phi(K), K)$ and $\sum^{N-1}_{k=0}|K_{t^N_{k+1}}-K_{t^N_k}|^2\in L^1_G(\Omega_T)$, we should assume here that $K_T\in L^2_G(\Omega_T)$.
\end{remark}
\begin{proposition}\label{boudiff}
Let $X$ be a q.s. continuous $G$-It\^o process such that 
\begin{equation}\label{gito}
X_t= x+\int^t_0 f_sds+\int^t_0h_sd\langle
B\rangle_s+\int^t_0g_sdB_s,\ 0\leq t\leq T,
\end{equation}
where $f$, $h$ and $g$ are elements in $M^p_G([0, T])$, $p>2$. Let $K\in M_I([0, T])\cap M^q_G([0, T])$ and $K_T\in L^q_G(\Omega_T)$, where $1/p+1/q=1$. Then, $\int^T_0X_tdK_t$ is an element in $L^1_G(\Omega_T)$. 
\end{proposition}
\noindent\textbf{Proof:} Given a sequence of refining partitions $\{\pi^N_{[0,T]}\}_{N\in\mathbb{N}}$, we construct sequence (\ref{constr}). By the definitions of stochastic integrals and the BDG type inequalities, one can verify that for each $t\in [0, T]$, $X_t\in L^p_G(\Omega_t)$. Therefore, for all $N\in\mathbb{N}$, $\mathcal{V}_{[0,T]}^N(X, K)\in L^1_G(\Omega_T)$. Applying the BDG type inequalities, we have
\begin{align*}
\bar{\mathbb{E}}[\sup_{s\leq u\leq t}|X_u-X_s|^p]
\leq C\bigg(|t-s|^{p-1}\bigg(\int^t_s(\bar{\mathbb{E}}[|f_u|^p]
&+\bar{\mathbb{E}}[|h_u|^p])du\bigg)\\
&+|t-s|^{\frac{p}{2}-1}\int^t_s\bar{\mathbb{E}}[|g_u|^p]du\bigg).
\end{align*}
Thus, 
\begin{align}\label{za}
&\ \bar{\mathbb{E}}[\sup_{k\in[0, N)\cap\mathbb{N}}\sup_{t_k\leq t\leq t_{k+1}}|X_t-X_{t^N_k}|^p]
\leq\bar{\mathbb{E}}\bigg[\sum_{k=0}^{N-1}\sup_{t^N_k\leq t\leq t^N_{k+1}}|X_t-X_{t^N_k}|^p\bigg]\notag\\
\leq&\ C\sum_{k=0}^{N-1}\bigg(\int^{t^N_{k+1}}_{t^N_k}(|t^N_{k+1}-t^N_k|^{p-1}(\bar{\mathbb{E}}[|f_t|^p]+\bar{\mathbb{E}}[|h_t|^p])
+|t^N_{k+1}-t^N_k|^{\frac{p}{2}-1}\bar{\mathbb{E}}[|g_t|^p])dt\bigg)\\
\leq&\ C\bigg(\mu(\pi^N_{[0,T]})^{p-1}\int^T_0(\bar{\mathbb{E}}[|f_t|^p]+\bar{\mathbb{E}}[|h_t|^p])dt+\mu(\pi^N_{[0,T]})^{\frac{p}{2}-1}\int^T_0\bar{\mathbb{E}}[|g_t|^p]dt\bigg).\notag
\end{align}
From the integrability of $f$, $h$ and $g$, we have
\begin{equation}\bar{\mathbb{E}}[\sup_{k\in[0, N)\cap\mathbb{N}}\sup_{t_k\leq t\leq t_{k+1}}|X_t-X_{t^N_k}|^p]
\leq CM(\mu(\pi^N_{[0,T]})^{p-1}+\mu(\pi^N_{[0,T]})^{\frac{p}{2}-1}).\notag
\end{equation}
For each $N\in\mathbb{N}$, we calculate
\begin{align*}
&\ \bigg|\mathcal{V}_{[0,T]}^N(X, K)-\int^T_0 X_tdK_t\bigg|\leq 
\int^T_0 \bigg|\sum^{N-1}_{k=0}X_{t^N_k}{\mathbb{1}}_{[t^N_k, t^N_{k+1})}(t)-X_t\bigg|dK_t\\
\leq&\ \sup_{0\leq t\leq T}\bigg|\sum^{N-1}_{k=0}X_{t^N_k}{\mathbb{1}}_{[t^N_k, t^N_{k+1})}(t)-X_t\bigg|K_T\leq K_T\sup_{k\in[0, N)\cap\mathbb{N}}\sup_{t_k\leq t< t_{k+1}}|X_t-X_{t^N_k}|.
\end{align*}
Consequently,
\[
\begin{array}{r@{\ }l}
\bar{\mathbb{E}}[|\mathcal{V}_{[0,T]}^N(X, K)-\int^T_0 X_tdK_t|]
\leq&\bar{\mathbb{E}}[K_T\sup\limits_{k\in[0, N)\cap\mathbb{N}}\sup\limits_{t_k\leq t< t_{k+1}}|X_t-X_{t^N_k}|]\\[6pt]
\leq&
(\bar{\mathbb{E}}[\sup\limits_{k\in[0, N)\cap\mathbb{N}}\sup\limits_{t_k\leq t< t_{k+1}}|X_t-X_{t^N_k}|^p])^{\frac{1}{p}}(\bar{\mathbb{E}}[K^q_T])^{\frac{1}{q}} \\[6pt]
\leq&CM(\mu(\pi^N_{[0,T]})^{p-1}+\mu(\pi^N_{[0,T]})^{\frac{p}{2}-1})^{\frac{1}{p}}\rightarrow 0,\  \textnormal{as}\ N\rightarrow +\infty.
\end{array}
\]
The desired result follows.\hfill{}$\square$
\subsection{An extension of $G$-It\^o's formula} 
\noindent For each $0\leq s\leq t\leq T$, consider a sum of a $G$-It\^o process and an increasing process $K$:
\[
X_t= X_s + \int^t_s f_udu+\int^t_sh_ud\langle
B\rangle_u+\int^t_sg_udB_u+K_t-K_s.
\]
\begin{lemma}\label{itolem} Let $\Phi\in \mathcal{C}^2(\mathbb{R})$ be a
real function with bounded and Lipschitz derivatives. Let $f$, $h$ and $g$ be bounded
processes in $M^2_G([0, T])$ and $K\in M_I([0, T])\cap M^2_G([0,T])$ satisfy for each $t\in[0, T]$, 
\begin{equation}\label{kconti}
\lim_{s\rightarrow t}\bar{\mathbb{E}}[|K_t-K_s|^2]=0.
\end{equation}
Then, 
\begin{align}\label{ITO}
\Phi(X_t)-\Phi(X_s)&=\int^t_s\frac{d\Phi}{dx}(X_u)f_udu+\int^t_s\frac{d\Phi}{dx}(X_u)h_ud\langle
B\rangle_u\notag\\
&+\int^t_s\frac{d\Phi}{dx}(X_u)g_udB_u+\int^t_s\frac{d\Phi}{dx}(X_u)dK_u\\
&+\frac{1}{2}\int^t_s\frac{d^2\Phi}{dx^2}(X_u)g_u^2d\langle
B\rangle_u,\ q.s..\notag
\end{align}
\end{lemma}
\noindent The proof of this lemma is based on previous results in Peng \cite{P3} (cf. Lemma 6.1 and Proposition 6.3 in Chapter III). To avoid redundancy, we first prove a reduced  lemma when $f=h=g\equiv0$ to show how the increasing process $K$ plays a role in this dynamic and then give a sketch to indicate some key points to combine the simple lemma with the previous results in Peng \cite{P3}.  
\begin{lemma}\label{reduced}
Let $\Phi\in \mathcal{C}^2(\mathbb{R})$ be a
real function with bounded and Lipschitz derivatives and $K\in M_I([0, T])\cap M^2_G([0,T])$. Then, 
\begin{align*}
\Phi(K_t)-\Phi(K_s)=
\int^t_s\frac{d\Phi}{dx}(K_u)dK_u,\ q.s..
\end{align*}
\end{lemma}
\noindent\textbf{Proof: } Consider a sequence of refining partitions $\{\pi^N_{[s,t]}\}_{N\in\mathbb{N}}$. For each $N\in\mathbb{N}$, from the second order Taylor expansion, we have
\begin{align*}
\Phi(K_t)-\Phi(K_s)&=\sum^{{N}-1}_{k=0}(\Phi(K_{t^{N}_{k+1}})-\Phi(K_{t^{N}_k}))\\
&=\sum^{{N}-1}_{k=0}\frac{d\Phi}{dx}(K_{t^{N}_{k}})(K_{t^{N}_{k+1}}-K_{t^{N}_k})
+\frac{1}{2}\sum^{{N}-1}_{k=0}\frac{d^2\Phi}{dx^2}(\xi^{N}_{k})(K_{t^{N}_{k+1}}-K_{t^{N}_k})^2,
\end{align*}
where $\xi^{N}_k$ satisfies $K_{t^{N}_k}\leq \xi^{N}_k\leq
K_{t^{N}_{k+1}}$, q.s.. For the first part, similar to that in Remark \ref{sequen}, we obtain 
\[
\lim_{N\rightarrow +\infty}\bigg|\sum^{{N}-1}_{k=0}\frac{d\Phi}{dx}(K_{t^{N}_{k}})(K_{t^{N}_{k+1}}-K_{t^{N}_k})-
\int^t_s\frac{d\Phi}{dx}(K_u)dK_u\bigg|=0,\ q.s..
\]
For the second part, because $\frac{d^2\Phi}{dx^2}$ is bounded and the quadratic variation of $K$ on $[0, T]$ is q.s. 0, then,  
\[\frac{1}{2}\sum^{{N}-1}_{k=0}\frac{d^2\Phi}{dx^2}(\xi^{N}_{k})(K_{t^{N}_{k+1}}-K_{t^{N}_k})^2\leq \frac{1}{2}M\sum^{{N}-1}_{k=0}(K_{t^{N}_{k+1}}-K_{t^{N}_k})^2\rightarrow 0,\ q.s.,\ \textnormal{as}\ N\rightarrow+\infty.
\]
The proof is complete. \hfill{}$\square$\\[6pt]
\noindent\textbf{Sketch of the proof of Lemma \ref{itolem}:}
To combine the result above with the ones in Peng \cite{P3}, we decompose $X$ into $M^X+K$, where $M^X$ denotes the $G$-It\^o part of $X$. Given a sequence of refining partitions $\{\pi^{2^N}_{[s,t]}\}_{N\in\mathbb{N}}$: for each $N\in\mathbb{N}$,
\[\pi^{2^N}_{[s,t]}=\{t^{2^N}_0,t^{2^N}_1 \ldots,
t^{2^N}_{2^N}\}=\{s,s+\delta,\ldots,s+{2^N}\delta=t\},\]
we have from the second order Taylor expansion
\begin{align*}
\Phi(X_t)-\Phi(X_s)&=\sum^{2^N-1}_{k=0}(\Phi(X_{t^{2^N}_{k+1}})-\Phi(X_{t^{2^N}_{k+1}}))\\
&=\sum^{{2^N}-1}_{k=0}\frac{d\Phi}{dx}(X_{t^{2^N}_{k}})(M^X_{t^{2^N}_{k+1}}-M^X_{t^{2^N}_k})
+\frac{1}{2}\sum^{2^{N}-1}_{k=0}\frac{d^2\Phi}{dx^2}(X_{t^{2^N}_{k}})(M^X_{t^{2^N}_{k+1}}-M^X_{t^{2^N}_k})^2\\
&+\sum^{2^{N}-1}_{k=0}\frac{d^2\Phi}{dx^2}(\xi^{2^N}_{k})(M^X_{t^{2^N}_{k+1}}-M^X_{t^{2^N}_k})(K_{t^{2^N}_{k+1}}-K_{t^{2^N}_k})
+\frac{1}{2}\sum^{2^{N}-1}_{k=0}\frac{d^2\Phi}{dx^2}(\xi^{2^N}_{k})(K_{t^{2^N}_{k+1}}-K_{t^{2^N}_k})^2\\
&+\frac{1}{2}\sum^{2^{N}-1}_{k=0}\bigg(\frac{d^2\Phi}{dx^2}(\xi^{2^N}_{k})-\frac{d^2\Phi}{dx^2}(X_{t^{2^N}_{k}})\bigg)(M^X_{t^{2^N}_{k+1}}-M^X_{t^{2^N}_k})^2
+\sum^{{2^N}-1}_{k=0}\frac{d\Phi}{dx}(X_{t^{2^N}_{k}})(K_{t^{2^N}_{k+1}}-K_{t^{2^N}_k})\\&=I^N_1+I^N_2+I^N_3+I^N_4+I^N_5+I^N_6,
\end{align*}
where $\xi^{2^N}_k$ satisfies $X_{t^{2^N}_k}\wedge X_{t^{2^N}_{k+1}} \leq \xi^{2^N}_k
\leq X_{t^{2^N}_k}\vee X_{t^{2^N}_{k+1}}$ q.s.. \\[6pt]
A key point in the proof is to verify the following convergences in $M^2_G([0, T])$:
\begin{align}\label{conver1}
\sum^{{2^N}-1}_{k=0}\frac{d\Phi}{dx}(X_{t^{2^N}_{k}}){\mathbb{1}}_{[t^{2^N}_k, t^{2^N}_{k+1})}(\cdot)\rightarrow \frac{d\Phi}{dx}(X_\cdot),\ \textnormal{as}\ N\rightarrow +\infty;
\end{align}
and
\begin{equation}\label{conver2}
\sum^{{2^N}-1}_{k=0}\frac{d^2\Phi}{dx^2}(X_{t^{2^N}_{k}}){\mathbb{1}}_{[t^{2^N}_k, t^{2^N}_{k+1})}(\cdot)\rightarrow\frac{d^2\Phi}{dx^2}(X_\cdot),\ \textnormal{as}\ N\rightarrow +\infty.
\end{equation}
For the $G$-It\^o part $M^X$, we deduce by the BDG type inequalities
\begin{align}\label{conver3}
\int^t_s\bar{\mathbb{E}}\bigg[\bigg|\sum^{{2^N}-1}_{k=0}M^X_{t^{2^N}_k}{\mathbb{1}}_{[t^{2^N}_k, t^{2^N}_{k+1})}(u)-M^X\bigg|^2\bigg]du\leq M|t-s|(\delta+\delta^2)\rightarrow 0,\ \textnormal{as}\ N\rightarrow +\infty.
\end{align}
For the increasing process $K$, thanks to assumption (\ref{kconti}), for each $u\in [s, t]$,
\begin{equation}\label{hell}
\lim_{N\rightarrow +\infty}\bar{\mathbb{E}}\bigg[\bigg|\sum^{{2^N}-1}_{k=0}K_{t^{2^N}_k}{\mathbb{1}}_{[t^{2^N}_k, t^{2^N}_{k+1})}(u)-K_u\bigg|^2\bigg]= 0.
\end{equation}
Moreover, 
\[
\int^t_s\bar{\mathbb{E}}\bigg[\bigg|\sum^{{2^N}-1}_{k=0}K_{t^{2^N}_k}{\mathbb{1}}_{[t^{2^N}_k, t^{2^N}_{k+1})}(u)\bigg|^2\bigg]du\leq\int^t_s\bar{\mathbb{E}}[K^2_u]du<+\infty.
\]
By Lebesgue's dominated convergence theorem to the integral on $[s, t]$, we deduce
\begin{equation}\label{conver4}
\lim_{N\rightarrow +\infty}\int^t_s\bar{\mathbb{E}}\bigg[\bigg|\sum^{{2^N}-1}_{k=0}K_{t^{2^N}_k}{\mathbb{1}}_{[t^{2^N}_k, t^{2^N}_{k+1})}(u)-K_u\bigg|^2\bigg]du= 0.
\end{equation}
Combining (\ref{conver3}) and (\ref{conver4}), (\ref{conver1}) and (\ref{conver2}) are readily obtained by the Lipschitz continuity of $\frac{d\Phi}{dx}$ and $\frac{d^2\Phi}{dx^2}$. Then, we can proceed similarly to Peng \cite{P3} to treat with $I^N_1$ and $I^N_2$. \\[6pt]
On the other hand, due to the boundedness of $\frac{d^2\Phi}{dx^2}$ and the boundedness and uniform continuity of paths $M^X_\cdot(\omega)$ and $K_\cdot(\omega)$ on $[0, T]$, for each $\omega\in A^c$, we can easily obtain that $I^N_3$ and $I^N_4$ q.s. vanish.\\[6pt]
For $I^N_5$, we calculate
\begin{align*}
|I^N_5|&\leq \frac{C}{2}\sum^{2^{N}-1}_{k=0}|\xi^{2^N}_{k}-X_{t^{2^N}_{k}}||M^X_{t^{2^N}_{k+1}}-M^X_{t^{2^N}_k}|^2\\
&\leq \frac{C}{2}\bigg(\sum^{2^{N}-1}_{k=0}|(\xi^1)^{2^N}_{k}-M^X_{t^{2^N}_{k}}||M^X_{t^{2^N}_{k+1}}-M^X_{t^{2^N}_k}|^2+\sum^{2^{N}-1}_{k=0}|(\xi^2)^{2^N}_{k}-K_{t^{2^N}_{k}}||M^X_{t^{2^N}_{k+1}}-M^X_{t^{2^N}_k}|^2\bigg),
\end{align*}
where $(\xi^1)^{2^N}_k$ satisfies $M^X_{t^{2^N}_k}\wedge M^X_{t^{2^N}_{k+1}} \leq (\xi^1)^{2^N}_k
\leq M^X_{t^{2^N}_k}\vee M^X_{t^{2^N}_{k+1}}$ and $(\xi^2)^{2^N}_k$ satisfies $K_{t^{2^N}_k}\leq (\xi^2)^{2^N}_k
\leq K_{t^{2^N}_{k+1}}$, q.s.. The result in Peng \cite{P3} shows that the first part converges to 0 in $M^2_G([0, T])$, whereas the second part vanishes as a result of the uniform continuity of paths $K_\cdot(\omega)$ on $[0, T]$, for all $\omega\in A^c$ and the q.s. boundedness of the quadratic variation of the $G$-It\^o part $M^X$.\\[6pt]
\noindent For $I^N_6$, it converges to $\int^t_s\frac{d\Phi}{dx}(X_u)dK_u$, q.s. by Definition \ref{increa}.\hfill{}$\square$
\begin{remark}
In the proof of the classical It\^o's formula, (\ref{conver1}) and (\ref{conver2}) can be verified directly by the pathwise continuity of $X$ and Lebesgue's dominated convergence theorem on the product space $[s, t]\times\Omega$. But in the $G$-framework, we lack such a theorem. In general, given an $X\in M^2_G([0, T])$, the sequence of step processes $$\bigg\{\sum^{{2^N}-1}_{k=0}X_{t^{2^N}_k}{\mathbb{1}}_{[t^{2^N}_k, t^{2^N}_{k+1})}(\cdot)\bigg\}_{N\in\mathbb{N}}$$ could not converge to $X$ in the sense of $M^2_G([0, T])$. Thus, (\ref{kconti}) is needed to ensure that (\ref{hell}) holds true.\\[6pt]
\noindent In fact, the left-hand side of (\ref{ITO}), particularly the term $\int^t_s \frac{d\Phi}{dx}(X_u)dK_u$, still belongs to $L^2_G(\Omega_t)$. 
A sufficient condition of this result is that $K_t\in L^2_G(\Omega_t)$, which can be verified by choosing a sequence such that $t_n\rightarrow t$ and for all $n\in \mathbb{N}$, $X_{t_n}\in L^2_G(\Omega_{t_n})$ (Remark \ref{remconna} ensures the existence of this sequence) and by deduction from assumption (\ref{kconti}).
\end{remark}
\noindent Similar to Theorem 6.5 of Peng \cite{P3}, we can extend $G$-It\^o's formula in Lemma \ref{itolem} to those $\Phi$ whose second derivatives $\frac{d^2\Phi}{dx^2}$ have polynomial growth. Unfortunately, this extension is at the cost of more restrictions on the increasing process $K$. 
\begin{theorem}\label{BT}
Let $\Phi\in \mathcal{C}^2(\mathbb{R})$ be a real function such that $\frac{d^2\Phi}{dx^2}$ satisfies the polynomial growth condition. Let $f$, $h$ and $g$ be bounded
processes in $M^2_G([0, T])$ and $K\in M_I([0, T])\cap M^2_G([0,T])$ 
satisfies that for each $t\in[0, T]$,
\begin{align*}
\lim_{s\rightarrow t}\bar{\mathbb{E}}[|K_t-K_s|^2]=0;
\end{align*} and
for any $p>2$, $\bar{\mathbb{E}}[K^p_T]< +\infty$.
Then,
\begin{align}\label{TO}
\Phi(X_t)-\Phi(X_s)&=\int^t_s\frac{d\Phi}{dx}(X_u)f_udu+\int^t_s\frac{d\Phi}{dx}(X_u)h_ud\langle
B\rangle_u\notag\\
&+\int^t_s\frac{d\Phi}{dx}(X_u)g_udB_u+\int^t_s\frac{d\Phi}{dx}(X_u)dK_u\\
&+\frac{1}{2}\int^t_s\frac{d^2\Phi}{dx^2}(X_u)g_u^2d\langle
B\rangle_u,\ q.s..\notag
\end{align}
\end{theorem}
\noindent\textbf{Proof:} By the same argument in the proof of Theorem 6.5 of Peng \cite{P3}, we can choose a sequence of functions $\Phi^N\in \mathcal{C}^2_0(\mathbb{R})$, such that for each $x\in\mathbb{R}$,
\begin{equation}\label{sb1}
|\Phi^N(x)-\Phi(x)|+\bigg|\frac{d\Phi^N}{dx}(x)-\frac{d\Phi}{dx}(x)\bigg|+\bigg|\frac{d^2\Phi^N}{dx^2}(x)-\frac{d^2\Phi}{dx^2}(x)\bigg|\leq \frac{C}{N}(1+|x|^k),
\end{equation}
where $C$ and $k$ are positive constants independent of $N$. Obviously, $\Phi^N$ satisfies the conditions in Lemma \ref{itolem}. Therefore, 
\begin{align}\label{appito}
\Phi^N(X_t)-\Phi^N(X_s)&=\int^t_s\frac{d\Phi^N}{dx}(X_u)f_udu+\int^t_s\frac{d\Phi^N}{dx}(X_u)h_ud\langle B\rangle_u\notag\\
&+\int^t_s\frac{d\Phi^N}{dx}(X_u)g_udB_u+\int^t_s\frac{d\Phi^N}{dx}(X_u)dK_u\\
&+\frac{1}{2}\int^t_s\frac{d^2\Phi^N}{dx^2}(X_u)g^2_ud\langle B\rangle_u.\notag
\end{align}
Borrowing the notation in the proof of Lemma \ref{itolem} and using the BDG type inequalities, we have 
\begin{equation}\label{sb2}
\bar{\mathbb{E}}[\sup_{0\leq t\leq T}|X_t|^{2k}]\leq C(\bar{\mathbb{E}}[\sup_{0\leq t\leq T}|M^X_t|^{2k}]+\bar{\mathbb{E}}[|K_T|^{2k}])< +\infty.
\end{equation}
Then,  from (\ref{sb1}) and (\ref{sb2}), we deduce that as $N\rightarrow+\infty$,
\begin{align}\label{troisconver}
\Phi^N(X_t)&\rightarrow\Phi(X_t),\ {\rm in}\ L^2_G(\Omega_t);\notag\\
\frac{d\Phi^N}{dx}(X_\cdot)&\rightarrow\frac{d\Phi}{dx}(X_\cdot),\ {\rm in}\ M^2_G([0, T]);\\
\frac{d^2\Phi^N}{dx^2}(X_\cdot)&\rightarrow\frac{d^2\Phi}{dx^2}(X_\cdot),\ {\rm in}\ M^2_G([0, T]).\notag
\end{align}
We can proceed as in Peng \cite{P3} to show that the terms on the right-hand side of (\ref{appito}), except $\int^t_s\frac{d\Phi^N}{dx}(X_u)dK_u$, converge to their corresponding terms in (\ref{TO}). To complete the proof, it suffices to show that for each $\omega\in A^c$, 
\begin{align*}
&\ \bigg|\int^t_s\frac{d\Phi^N}{dx}(X_u(\omega))dK_u(\omega)-\int^t_s\frac{d\Phi}{dx}(X_u(\omega))dK_u(\omega)\bigg|\\
\leq&\ \frac{C}{N}\int^t_s (1+|X_u(\omega)|^k)dK_u(\omega)\leq \frac{C}{N}(1+M^k_{\omega})K_T(\omega)\rightarrow 0,\ \textnormal{as}\ N\rightarrow +\infty,
\end{align*}
by the continuity and boundedness of paths $X_\cdot(\omega)$ and $K_\cdot(\omega)$ on $[0, T]$. \hfill{}$\square$
\begin{remark}
If $|\frac{d^2\Phi}{dx^2}(x)|\leq C(1+|x|^{k})$, for some $k\geq 1$, then the condition on $K$ could be weakened to $\bar{\mathbb{E}}[|K_T|^{2(k+3)}]<+\infty$. 
\end{remark}
\begin{remark}
Following exactly the procedure above, we can obtain a similar result when a bounded variation process $K_1-K_2$ appears in the dynamic.
\end{remark}
\section{Reflected $G$-Brownian motion}
\noindent Before moving to the main result of this paper, we first consider a reduced RGSDE, that is, taking $f=h\equiv0$ and $g\equiv1$, only a $G$-Brownian motion and an increasing process drive the dynamic on the right-hand side of (\ref{ABC}). In what follows, we establish the solvability to the RGSDE of this type, i.e., the existence and uniqueness of reflected $G$-Brownian Motion.\\[6pt]
\noindent Let $y$ be a real valued continuous function on $[0, T]$ with $y_0 \geq 0$. It is well-known that there exists a unique pair $(x, k)$ of functions on $[0, T]$ such that
$x=y+k$, where $x$ is positive, $k$ is an increasing and continuous function that starts from 0. Moreover, the Riemann-Stieltjes integral $\int^T_0 x_t dk_t = 0$. The solution to this Skorokhod problem on $[0, T]$ is given by
\begin{equation}\label{ssp}
\left\{\begin{aligned} 
x_t &= y_t+k_t;\\
k_t &= \sup_{s\leq t} x^-_s,
\end{aligned}\right.
\end{equation}
which is explicit and unique. 
\begin{theorem}\label{RBT}
For any $p\geq 1$, there exists a unique pair of processes $(X, K)$ in $M^p_G([0, T])\times (M_I([0, T])\cap M^p_G([0, T]))$, such that
\begin{equation}\label{RB}
X_t=B_t+K_t,\ 0\leq t\leq T,\ \ q.s.,
\end{equation}
where (a) $K_0=0$; (b) $X$ is positive; and (c) $\int^T_0 X_t dK_t=0$, q.s..
\end{theorem}
\noindent\textbf{Proof:} With the help of (\ref{ssp}), we define a pair of processes $(X, K)$ pathwisely on $[0, T]$: 
\begin{equation}\label{RBS}
\left\{\begin{aligned} 
X_t(\omega) &= B_t(\omega)+K_t(\omega);\\
K_t(\omega) &= \sup\limits_{s\leq t} B^-_s(\omega).
\end{aligned}\right.
\end{equation}
Obviously, $K\in M_I([0, T])$ and (a), (b) and (c) are satisfied. Therefore, to complete the proof, we need only verify that $K\in M^p_G([0,T])$.\\[6pt]
Because for all $1\leq p'<p$, $M^{p'}_G([0, T])\subset M^{p}_G([0, T])$, we can assume that $p>2$ without loss of generality.  Given a sequence of partitions $\{\pi^N_{[0,T]}\}_{N\in\mathbb{N}}$,
we set
\begin{align*}
(B^-_t)^N(\omega):=\sum^{N-1}_{k=0} B^-_{t^N_k}
(\omega){\mathbb{1}}_{[t^N_k,t^N_{k+1})}(t),\ 0\leq t\leq T;
\end{align*}
and
\[
\sup\limits_{0\leq s\leq t} (B^-_s)^N := \sum^{N-1}_{k=0} \max_{l\in\{0, 1,\ldots, k\}} B^-_{t^N_l}
{\mathbb{1}}_{[t^N_k,t^N_{k+1})}(t),\ 0\leq t\leq T.
\]
We observe that both $((B^-_t)^N)_{0\leq t\leq T}$ and $(\sup\limits_{0\leq s\leq t} (B^-_s)^N)_{0\leq t\leq T}$ are step processes in $M^p_G([0,T])$. Because 
\begin{align*}
&\ \bar{\mathbb{E}}[|\sup_{0\leq s\leq t}(B^-_s)^N-\sup_{0\leq s\leq t}B^-_s|^{p}]
\leq\bar{\mathbb{E}}[\sup_{0\leq s\leq t}|(B^{-}_s)^N-B^-_s|^{p}]\\
\leq&\ \bar{\mathbb{E}}[\sup_{0\leq t\leq T}|B_t^N-B_t|^{p}]
\leq\bar{\mathbb{E}}[\sup_{k\in\mathbb{N}\cap[0, N)}\sup_{t_k\leq t< t_{k+1}}|B_t-B_{t^N_k}|^p],
\end{align*}
then, letting $f=h\equiv0$ and $g\equiv1$ in (\ref{za}), we obtain 
\begin{align*}
\bar{\mathbb{E}}[|\sup_{0\leq s\leq t}(B^-_s)^N-\sup_{0\leq s\leq t}B^-_s|^p]\leq C \mu(\pi^N_{[0,T]})^{\frac{p}{2}-1}
\rightarrow 0,\ \textnormal{as}\ N\rightarrow+\infty,
\end{align*}
which shows that $(\sup\limits_{0\leq s\leq t} (B^-_s)^N)_{0\leq t\leq T}$ converges to $K$ in $M^p_G([0, T])$. \\[6pt]
On the other hand, the uniqueness of such a pair $(X, K)$ is inherited from the solution to the Skorokhod problem pathwisely. The proof is complete.\hfill$\square$
\begin{remark}
We call the process $X$ in Theorem \ref{RBT} a $G$-reflected Brownian motion on the half-line $[0, +\infty)$.
\end{remark}
\noindent Furthermore, if the $G$-Brownian motion $B$ is replaced by some $G$-It\^o process, we have the following statement similar to Theorem \ref{RBT}.
\begin{theorem}\label{BXT}
For some $p>2$, consider a q.s. continuous $G$-It\^o process $Y$ defined in the form of (\ref{gito}) whose coefficients are all elements in $M^p_G([0,T])$.
Then, there exists a unique pair of processes $(X, K)$ in $M^{p}_G([0, T])\times (M_I([0, T])\cap M^{p}_G([0, T]))$ such that
\begin{equation}\label{suchthat}
X_t=Y_t+K_t,\ 0\leq t\leq T,\ \ q.s.,
\end{equation}
where (a) $X$ is positive; (b) $K_0=0$; and (c) $\int^T_0 X_t dK_t=0$, q.s..
\end{theorem}
\noindent  We omit the proof, because it is an analogue to the proof above and deduced mainly by the integrability of the coefficients of $Y$ and (\ref{za}).
\section{Scalar valued RGSDEs}
\noindent We state our main result in this section by giving the existence and uniqueness of the solutions to the scalar valued RGSDEs with Lipschitz coefficients. Additionally, a comparison theorem is given at the end of this paper.
\subsection{Formulation to RGSDEs}
\noindent We consider the following scalar valued RGSDE:
\begin{equation}\label{1}
X_t=x+\int^t_0 f_s(X_s)ds+\int^t_0 h_s(X_s)d\langle B \rangle_s +\int^t_0 g_s(X_s) dB_s +K_t,\ 0\leq t\leq T,\ q.s.,
\end{equation}
where\\[-6pt]
\begin{mlist}
\item The initial condition $x\in \mathbb{R}$;
\item For some $p>2$, the coefficients $f$, $h$ and $g: \Omega\times[0,T]\times\mathbb{R}\rightarrow\mathbb{R}$ are given functions that satisfy for each $x\in\mathbb{R}$, $f_\cdot(x)$, $h_\cdot(x)$ and $g_\cdot(x)\in M^p_G([0,T])$;
\item The coefficients $f$, $h$ and $g$ that satisfy a Lipschitz condition, i.e., for each $t\in[0, T]$ and $x, x'\in\mathbb{R}$, $|f_t(x)-f_t(x')|+|h_t(x)-h_t(x')|+|g_t(x)-g_t(x')|\leq C_L|x-x'|$, q.s.;
\item The obstacle is a $G$-It\^o process whose coefficients are all elements in $M^p_G([0,T])$,
and we shall always assume that $S_0\leq x$, q.s..
\end{mlist}
The solution to the RGSDE (\ref{1}) is a pair of processes $(X,K)$ that take values both in $\mathbb{R}$ and satisfy:\\[-6pt]
\begin{mlist}\addtocounter{nelist}{4} 
\item $X\in M^p_G([0,T])$ and $X_t\geq S_t$, $0\leq t\leq T$, q.s.;
\item $K\in M_I([0, T])\cap M^p_G([0,T])$ and $K_0=0$, q.s.;
\item $\int^T_0(X_t-S_t)dK_t=0$, q.s..
\end{mlist}
\subsection{Some a priori estimates and the uniqueness result}
\noindent Let $(X, K)$ be a solution to (\ref{1}). Replacing $Y_t$ by 
$x+\int^t_0 f_s(X_s)ds+\int^t_0 h_s(X_s)d\langle B \rangle_s +\int^t_0 g_s(X_s) dB_s-S_t$ and $X_t$ by $X_t-S_t$ in (\ref{suchthat}), we have the following representation of $K$:
\begin{align}\label{represen}
K_t=\sup_{0\leq s\leq t}\bigg(x+\int^s_0 f_u(X_u)du&+\int^s_0 h_u(X_u)d\langle B \rangle_u\\
&+\int^s_0 g_u(X_u) dB_u-S_s\bigg)^-,\ 0\leq t\leq T,\ q.s..\notag
\end{align}
\noindent We now give an a priori estimate on the uniform norm of the solution.
\begin{proposition}\label{abcde}
Let $(X, K)$ be a solution to (\ref{1}). Then, there exists a constant $C>0$ such that
\begin{align*}
\bar{\mathbb{E}}[\sup_{0\leq t\leq T}|X_t|^p]+\bar{\mathbb{E}}[K^p_T]
\leq C\bigg(|x|^p
&+\int^T_0(\bar{\mathbb{E}}[|f_t(0)|^p]\\
&+\bar{\mathbb{E}}[|h_t(0)|^p]
+\bar{\mathbb{E}}[|g_t(0)|^p])dt
+\bar{\mathbb{E}}[\sup_{0\leq t\leq T}|S^+_t|^p]\bigg).
\end{align*}
\end{proposition}
\noindent\textbf{Proof:} 
As $X$ is the solution to (\ref{1}), we obtain
\begin{equation}\label{est1}
\begin{array}{r@{\ }l}
\bar{\mathbb{E}}[\sup\limits_{0\leq s\leq t}|X_s|^p]
\leq&\bar{\mathbb{E}}[\sup\limits_{0\leq s\leq t}|x+\int^s_0 f_u(X_u)du+\int^s_0 h_u(X_u)d\langle B \rangle_u +\int^s_0 g_u(X_u) dB_u+K_s|^p]\\[6pt]
\leq& C(|x|^p+\bar{\mathbb{E}}[\sup\limits_{0\leq s\leq t}|\int^s_0 f_u(X_u)du|^p]
+\bar{\mathbb{E}}[\sup\limits_{0\leq s\leq t}|\int^s_0 h_u(X_u)d\langle B \rangle_u|^p]\\[6pt]
+&\bar{\mathbb{E}}[\sup\limits_{0\leq s\leq t}|\int^s_0 g_u(X_u)dB_u|^p]+\bar{\mathbb{E}}[|K_t|^p]).\end{array}
\end{equation}
In a similar way to (\ref{est1}), from the representation of $K$ (\ref{represen}), we have
\begin{equation}\label{est2}
\begin{array}{r@{\ }l}
\bar{\mathbb{E}}[K^p_t]\leq&\bar{\mathbb{E}}[\sup\limits_{0\leq s\leq t}(( x+\int^s_0 f_u(X_u)du+\int^s_0 h_u(X_u)d\langle B \rangle_u +\int^s_0 g_u(X_u) dB_u-S_s)^-)^p]\\[6pt]
\leq&\bar{\mathbb{E}}[\sup\limits_{0\leq s\leq t}(( x+\int^s_0 f_u(X_u)du+\int^s_0 h_u(X_u)d\langle B \rangle_u +\int^s_0 g_u(X_u) dB_u-S^+_s)^-)^p]\\[6pt]
\leq&\bar{\mathbb{E}}[\sup\limits_{0\leq s\leq t}|x+\int^s_0 f_u(X_u)du+\int^s_0 h_u(X_u)d\langle B \rangle_u +\int^s_0 g_u(X_u) dB_u-S^+_s|^p]\\[6pt]
\leq& C(|x|^p+\bar{\mathbb{E}}[\sup\limits_{0\leq s\leq t}|\int^s_0 f_u(X_u)du|^p]
+\bar{\mathbb{E}}[\sup\limits_{0\leq s\leq t}|\int^s_0 h_u(X_u)d\langle B \rangle_u|^p]\\[6pt]
+&\bar{\mathbb{E}}[\sup\limits_{0\leq s\leq t}|\int^s_0 g_u(X_u)dB_u|^p]+\bar{\mathbb{E}}[\sup\limits_{0\leq s\leq t}|S^+_s|^p]).
\end{array}
\end{equation}
Combining (\ref{est1}) and (\ref{est2}) and applying BDG type inequalities, we get 
\begin{equation*}
\begin{array}{r@{\ }l}
\bar{\mathbb{E}}[\sup\limits_{0\leq s\leq t}|X_s|^p]+\bar{\mathbb{E}}[K^p_t]\leq C(|x|^p+&\int^t_0 (\bar{\mathbb{E}}[|f_s(X_s)|^p]\\[6pt]
+&\bar{\mathbb{E}}[|h_s(X_s)|^p]+\bar{\mathbb{E}}[|g_s(X_s)|^p])ds+\bar{\mathbb{E}}[\sup\limits_{0\leq s\leq t}|S^+_s|^p].
\end{array}
\end{equation*}
By assumption (A3), we calculate
\begin{equation}\label{est3}
\begin{array}{r@{\ }l}
\bar{\mathbb{E}}[\sup\limits_{0\leq s\leq t}|X_s|^p]+\bar{\mathbb{E}}[K^p_t] 
\leq& C(|x|^p
+\int^t_0(\bar{\mathbb{E}}[(|f_s(0)|+C_L|X_s|)^p]
+\bar{\mathbb{E}}[(|h_s(0)|+C_L|X_s|)^p]\\[6pt]
+&\bar{\mathbb{E}}[(|g_s(0)+C_L|X_s|)^p])ds+\bar{\mathbb{E}}[\sup\limits_{0\leq s\leq t}|S^+_s|^p]\\[6pt]
\leq&C(|x|^p
+\int^t_0 (\bar{\mathbb{E}}[|f_s(0)|^p]+\bar{\mathbb{E}}[|h_s(0)|^p]
+\bar{\mathbb{E}}[|g_s(0)|^p])ds\\[6pt]
+&\bar{\mathbb{E}}[\sup\limits_{0\leq s\leq t}|S^+_s|^p]+\int^t_0\bar{\mathbb{E}}[|X_s|^p]ds)\\[6pt]
\leq&C(|x|^p
+\int^T_0 (\bar{\mathbb{E}}[|f_t(0)|^p]+\bar{\mathbb{E}}[|h_t(0)|^p]
+\bar{\mathbb{E}}[|g_t(0)|^p])dt\\[6pt]
+&\bar{\mathbb{E}}[\sup\limits_{0\leq t\leq T}|S^+_t|^p]+\int^t_0\bar{\mathbb{E}}[\sup\limits_{0\leq u\leq s}|X_u|^p]ds).
\end{array}
\end{equation}
Applying Gronwall's lemma to $\bar{\mathbb{E}}[\sup\limits_{0\leq s\leq t}|X_s|^p]$, we deduce
\begin{equation}\label{est4}
\begin{array}{r@{\ }l}
\bar{\mathbb{E}}[\sup\limits_{0\leq s\leq t}|X_t|^p]\leq C\bigg(|x|^p
+&\int^T_0(\bar{\mathbb{E}}[|f_t(0)|^p]
+\bar{\mathbb{E}}[|h_t(0)|^p]\\[6pt]
+&\bar{\mathbb{E}}[|g_t(0)|^p])dt
+\bar{\mathbb{E}}[\sup\limits_{0\leq t\leq T}|S^+_t|^p]\bigg),\ 0\leq t\leq T.
\end{array}
\end{equation}
Putting (\ref{est4}) into (\ref{est3}), the result follows.\hfill{}$\square$\\[6pt]
\noindent In the following theorem, we estimate the variation in the solutions induced by a variation in the coefficients and the obstacle processes.
\begin{theorem}\label{345}
Let $(x^1, f^1, h^1, g^1, S^1)$ and $(x^2, f^2, h^2, g^2, S^2)$ be two sets of coefficients that satisfy the assumptions (A1)-(A4) and $(X^i, K^i)$ the solution to the RGSDE corresponding to $(x^i, f^i, h^i, g^i, S^i)$, $i=1, 2$. Define
\begin{align*}
\Delta x:= x^1-x^2,\ \Delta f := f^1-f^2,\ &\Delta h := h^1-h^2,\ \Delta g := g^1-g^2;\\
\Delta S := S^1-S^2,\ \Delta X := X^1&-X^2,\ \Delta K := K^1-K^2.
\end{align*}
Then there exists a constant $C>0$ such that
\begin{align*}\bar{\mathbb{E}}[\sup_{0\leq t\leq T}|\Delta X_t|^p]
\leq C\bigg(|\Delta x|^p
+\int^T_0(\bar{\mathbb{E}}[|\Delta f_t(X_t^1)|^p]
+&\bar{\mathbb{E}}[|\Delta h_t(X_t^1)|^p]\\
+&\bar{\mathbb{E}}[|\Delta g_t(X_t^1)|^p])dt+\bar{\mathbb{E}}[\sup_{0\leq t\leq T}|\Delta S_t|^p]\bigg).
\end{align*}
\end{theorem}
\noindent\textbf{Proof:} Defining
\[
\begin{array}{r@{\ }l}
(M^X)^i_t:=x^i+\int^t_0 f^i_s(X^i_s)ds+\int^t_0 h^i_s(X^i_s)d\langle B \rangle_s +\int^t_0 g^i_s(X^i_s) dB_s,\ 0\leq t\leq T,\ i=1, 2;\\
\end{array}
\]
and
\[\Delta M^X:=(M^X)^1-(M^X)^2,\]
we calculate in a similar way to the proof of Proposition \ref{abcde}
\[
\begin{array}{r@{\ }l}
\bar{\mathbb{E}}[\sup\limits_{0\leq s\leq t}|(\Delta M^X)_s|^p]\leq& \bar{\mathbb{E}}[\sup\limits_{0\leq s\leq t}|\Delta x
+\int^s_0 (f^1_u(X^1_u)-f^2_u(X^2_u))du\\[6pt]
+&\int^s_0 (h^1_u(X^1_u)-h^2_u(X^2_u))d\langle B\rangle_u
+\int^s_0 (g^1_u(X^1_u)-g^2_u(X^2_u))dB_u|^p]\\[6pt]
\leq&\bar{\mathbb{E}}[\sup\limits_{0\leq s\leq t}|\Delta x
+\int^s_0 \Delta f_u(X_u^1) du
+\int^s_0 (f^2_u(X^1_u)-f^2_u(X^2_u))du\\[6pt]
+&\int^s_0 \Delta h_u(X_u^1)d\langle B\rangle_u
+\int^s_0 (h^2_u(X^1_u)-h^2_u(X^2_u))d\langle B\rangle_u\\[6pt]
+&\int^s_0 \Delta g_u(X_u^1)dB_u
+\int^s_0 (g^2_u(X^1_u)-g^2_u(X^2_u))dB_u|^p]\\[6pt]
\leq& C(|\Delta x|^p
+\int^t_0(\bar{\mathbb{E}}[|\Delta f_s(X_s^1)|^p]
+\bar{\mathbb{E}}[|\Delta h_s(X_s^1)|^p]\\[6pt]+&\bar{\mathbb{E}}[|\Delta g_s(X_s^1)|^p])ds
+\int^t_0\bar{\mathbb{E}}[|\Delta X_s|^p]ds)
\end{array}
\]
and
\begin{equation}\label{est6}
\begin{array}{r@{\ }l}
\bar{\mathbb{E}}[\sup\limits_{0\leq s\leq t}|\Delta K_s|^p]&=\bar{\mathbb{E}}[\sup\limits_{0\leq s\leq t}|\sup\limits_{0\leq u\leq s}((M^X)^1_u-S^1_u)^--\sup\limits_{0\leq u\leq s}((M^X)^2_u-S^2_u)^-|^p]\\[6pt]
&\leq\bar{\mathbb{E}}[\sup\limits_{0\leq s\leq t}|\sup\limits_{0\leq u\leq s}|{((M^X)^1_u-S^1_u)^-}-((M^X)^2_u-S^2_u)^-||^p]\\[6pt]
&=\bar{\mathbb{E}}[\sup\limits_{0\leq s\leq t}|{((M^X)^1_s-S^1_s)^-}-((M^X)^2_s-S^2_s)^-|^p]\\[6pt]
&\leq\bar{\mathbb{E}}[\sup\limits_{0\leq s\leq t}|((M^X)^1_s-S^1_s)-((M^X)^2_s-S^2_s)|^p]\\[6pt]
&\leq C(\bar{\mathbb{E}}[\sup\limits_{0\leq s\leq t}|\Delta (M^X)_s|^p]+\bar{\mathbb{E}}[\sup\limits_{0\leq s\leq t}|\Delta S_s|^p]).
\end{array}
\end{equation}
Then, we have
\[
\begin{array}{r@{\ }l}
\bar{\mathbb{E}}[\sup\limits_{0\leq s \leq t}|\Delta X_s|^p]
\leq& \bar{\mathbb{E}}[\sup\limits_{0\leq s\leq t}|(\Delta M^X)_s+\Delta K_s|^p|]\\[6pt]
\leq& C(\bar{\mathbb{E}}[\sup\limits_{0\leq s\leq t}|(\Delta M^X)_s|^p]+\bar{\mathbb{E}}[\sup\limits_{0\leq s\leq t}|\Delta K_s|^p])\\[6pt]
\leq& C(|\Delta x|^p
+\int^t_0(\bar{\mathbb{E}}[|\Delta f_s(X_s^1)|^p]
+\bar{\mathbb{E}}[|\Delta h_s(X_s^1)|^p]\\[6pt]+&\bar{\mathbb{E}}[|\Delta g_s(X_s^1)|^p])ds
+\bar{\mathbb{E}}[\sup\limits_{0\leq s\leq t}|\Delta S_s|^p]+\int^t_0\bar{\mathbb{E}}[|\Delta X_s|^p]ds).\\[6pt]
\leq& C(|\Delta x|^p
+\int^T_0(\bar{\mathbb{E}}[|\Delta f_t(X_t^1)|^p]
+\bar{\mathbb{E}}[|\Delta h_t(X_t^1)|^p]+\bar{\mathbb{E}}[|\Delta g_t(X_t^1)|^p])dt\\[6pt]
+&\bar{\mathbb{E}}[\sup\limits_{0\leq t\leq T}|\Delta S_t|^p]
+\int^t_0\bar{\mathbb{E}}[\sup\limits_{0\leq u\leq s}|\Delta X_u|^p]ds).\\[6pt]
\end{array}
\]
Gronwall's lemma gives the desired result.\hfill{}$\square$\\[6pt]
We deduce immediately the following uniqueness result by taking $x^1=x^2$, $f^1=f^2$, $h^1=h^2$, $g^1=g^2$ and $S^1=S^2$ in Theorem \ref{345}.
\begin{theorem}
Under assumptions (A1)-(A4), there exists at most one solution in $M^p_G([0,T])\times (M_I([0, T])\cap M^p_G([0,T]))$ to the RGSDE (\ref{1}).
\end{theorem}
\subsection{Existence result}
\noindent We now turn to the following existence result for the RGSDE (\ref{1}). The proof will be based on a Picard iteration.
\begin{theorem}\label{rue}
Under assumptions (A1)-(A4), there exists a unique solution in $M^p_G([0,T])\times (M_I([0, T])\cap M^p_G([0,T]))$ to the RGSDE (\ref{1}).
\end{theorem}
\noindent\textbf{Proof:}
We set $X^0=x$ and $K^0=0$. For each $n>0$, $X^{n+1}$ is given by recurrence:
\begin{align}\label{picard}
&X^{n+1}_t=x+\int^t_0 f_s(X^n_s)ds+\int^t_0 h_s(X^n_s)d\langle B \rangle_s +\int^t_0 g_s(X^n_s) dB_s + K^{n+1}_t,\ 0\leq t\leq T,
\end{align}
where\\[6pt]
$
\begin{array}{l}
(a)\ X^{n+1}\in M^p_G([0, T]),\ X^{n+1}_t\geq S_t,\ 0\leq t\leq T,\ q.s.;\\[6pt]
(b)\ K^{n+1}\in M_I([0, T])\cap M^p_G([0, T]),\ K^{n+1}_0=0,\ q.s.;\\[6pt]
(c)\ \int^T_0(X_t^{n+1}-S_t)dK^{n+1}_t=0,\ q.s..
\end{array}
$\\[6pt]
Substituting $X^{n+1}$ by $\tilde{X}^{n+1}+S_t$ on the left-hand side of (\ref{picard}), we know that $(X^{n+1}, K^{n+1})$ is well defined in $M^p_G([0,T])\times (M_I([0, T])\cap M^p_G([0,T]))$ by Theorem \ref{BXT}.\\[6pt]
First, we establish an a priori estimate uniform in $n$ for $\{\bar{\mathbb{E}}[\sup\limits_{0\leq t\leq T}|X^n_t|^p]\}_{n\in \mathbb{N}}$. In a similar way to (\ref{est3}), we have
\begin{align*}
\bar{\mathbb{E}}[\sup_{0\leq s\leq t}|X^{n+1}_t|^p]
\leq C\bigg(|x|^p
&+\int^T_0 (\bar{\mathbb{E}}[|f_t(0)|^p]+\bar{\mathbb{E}}[|h_t(0)|^p]\\
&+\bar{\mathbb{E}}[|g_t(0)|^p])dt
+\bar{\mathbb{E}}[\sup_{0\leq t\leq T}|S^+_t|^p]
+\int^t_0 \bar{\mathbb{E}}[\sup_{0\leq u\leq s}|X^n_u|^p]ds\bigg).
\end{align*}
By recurrence, it is easy to verify that for all $n\in\mathbb{N}$,
\[\bar{\mathbb{E}}[\sup_{0\leq s\leq t}|X^n_s|^p]\leq p(t),\ 0\leq t\leq T,\]
where $p(\cdot)$ is the solution  to the following ordinary differential equation:
\[p(t)=C\bigg(|x|^p
+\int^T_0 (\bar{\mathbb{E}}[|f_t(0)|^p]+\bar{\mathbb{E}}[|h_t(0)|^p]
+\bar{\mathbb{E}}[|g_t(0)|^p])dt
+\bar{\mathbb{E}}[\sup_{0\leq t\leq T}|S^+_t|^p]
+\int^t_0p(s)ds\bigg);\]
and $p(\cdot)$ is continuous and thus, bounded on $[0, T]$.\\[6pt]
Secondly, for each $n$ and $m\in \mathbb{N}$, we define
\[u^{n+1,m}_t:=\bar{\mathbb{E}}[\sup_{0\leq s \leq t}|X^{n+m+1}_s-X^{n+1}_s|^p],\ 0\leq t\leq T.\] 
Following the procedures in the proof of Theorem \ref{345}, we have
\[u^{n+1,m}_t \leq C\int^t_0u^{n,m}_sds.\]
Set
\[
v^n_t := \sup_{m\in\mathbb{N}} u^{n, m}_t,\ 0\leq t\leq T,
\]
then
\[0\leq u^{n+1, m}_t\leq C\sup_{m\in\mathbb{N}}\int^t_0u^{n, m}_sds\leq C\int^t_0\sup_{m\in\mathbb{N}}u^{n, m}_sds =C\int^t_0 v^n_s ds.\]
Taking the supremum over all $m\in\mathbb{N}$ on the left-hand side, we obtain
\[0\leq v^{n+1}_t=\sup_{m\in\mathbb{N}}u^{n+1,m}_t \leq C\int^t_0 v^n_s ds.\]
Finally, we define
\begin{equation}
\alpha_t := \limsup_{n\rightarrow+\infty} v^n_t,\ 0\leq t\leq T.\notag
\end{equation}
It is easy to find that $v^n_t \leq Cp(t)$, where $C$ is independent of $n$. By the Fatou-Lebesgue theorem,  we have
\[0\leq \alpha_t \leq C\int^t_0 \alpha_s ds.\]
Gronwall's lemma gives
\begin{equation}
\alpha_t=0,\ 0\leq t\leq T,\notag
\end{equation}
which implies that $\{X^n\}_{n\in\mathbb{N}}$ is a Cauchy sequence under the norm $(\mathbb{E}[\sup\limits_{0\leq t\leq T}|\cdot|^p])^{\frac{1}{p}}$, whose limit is certainly in $M^p_G([0, T])$. We denote the limit by $X$ and set
\[
K_t:=\sup_{0\leq s\leq t}\bigg( x+\int^s_0 f_u(X_u)du+\int^s_0 h_u(X_u)d\langle B \rangle_u +\int^s_0 g_u(X_u) dB_u-S_s\bigg)^-,\ 0\leq t\leq T.
\]
Obviously, the pair of processes $(X, K)$ satisfies (A5) - (A7).
We notice that
\begin{align*}
&\bar{\mathbb{E}}\bigg[\sup_{0\leq t\leq T}\bigg|\int^t_0(f_s(X^n_s)-f_s(X_s))ds\bigg|^p\bigg]\leq C\int^T_0\bar{\mathbb{E}}[|X^n_t-X_t|^p]dt;\\
&\bar{\mathbb{E}}\bigg[\sup_{0\leq t\leq T}\bigg|\int^t_0(h_s(X^n_s)-h_s(X_s))d\langle B\rangle_s\bigg|^p\bigg]\leq C\int^T_0\bar{\mathbb{E}}[|X^n_t-X_t|^p]dt;\\
&\bar{\mathbb{E}}\bigg[\sup_{0\leq t\leq T}\bigg|\int^t_0(g_s(X^n_s)-g_s(X_s))dB_s\bigg|^p\bigg]\leq C\int^T_0\bar{\mathbb{E}}[|X^n_t-X_t|^p]dt.
\end{align*}
Then, one can verify that $K^n$ converges to $K$ in $M^p_G([0, T])$ following the steps of (\ref{est6}). We conclude that the pair of processes $(X,K)$, well defined in $M^p_G([0, T])\times (M_I([0, T])\cap M^p_G([0, T]))$, is a solution to the RGSDE (\ref{1}).\hfill{}$\square$
\begin{remark}
Unlike a classical RSDE or RBSDE, the constraint process $S$ here is assumed to be a $G$-It\^o process instead of  a continuous process with $\bar{\mathbb{E}}[\sup\limits_{0\leq t\leq T}(S^+_t)^2]\leq +\infty$ (cf. El Karoui et al. \cite{EL}). In fact,  this is a sufficient condition to ensure that $K^{n+1}$ is still a $M^p_G([0, T])$ process in (\ref{picard}) by Theorem \ref{BXT}, which may be weakened to: 
\[\bar{\mathbb{E}}[\sup_{s\leq u\leq t}|S_u-S_s|^p]\leq C|t-s|^{\frac{p}{2}}.\] 
\end{remark}
\begin{remark}
We could also consider (\ref{1}) with less regularity assumptions on the coefficients $f$, $h$, $g$ and the obstacle $S$ under a family $\mathcal{P}_W$ of local martingale measures by using the approach introduced in Soner et al. \cite{STZ1, STZ3, STZ2}. The only problem in this case is the aggregation of the processes in the Picard iteration (\ref{picard}) for the proof of existence. Adapting to the assumptions of Theorem 2.2 in Nutz \cite{N}, we assume in addition that we work under Zermelo-Fraenkel set theory with the axiom of choice and the continuum hypothesis, then the stochastic integral of It\^o's type $\int^t_0g_s(X^n_s)dB_s$ can be well aggregated at each step of the recurrence. Thus, we can define from (\ref{ssp}) a universal pair $(X^{n+1}, K^{n+1})$ to make (\ref{picard}) $\mathbb{P}$-a.s. hold for all $\mathbb{P}\in\mathcal{P}_W$.  Following the argument in the proof of Theorem \ref{rue} under each $\mathbb{P}\in\mathcal{P}_W$, there exists a pair $(X^\mathbb{P}, K^\mathbb{P})$ such that (\ref{1}) holds true $\mathbb{P}$-a.s. and $(X^{n}, K^{n})$ converges to $(X^\mathbb{P}, K^\mathbb{P})$ in probability measure $\mathbb{P}$. By Lemma 2.5 in Nutz \cite{N}, there exists a universal pair $(X, K)$ such that $(X, K)=(X^\mathbb{P}, K^\mathbb{P})$, $\mathbb{P}$-a.s., which solves (\ref{1}) under each $\mathbb{P}\in\mathcal{P}$ and thus, in the (weaker) sense of $\mathcal{P}_W$-q.s.. 
\end{remark}
\begin{remark}
In contrast with the fact mentioned in Remark 3.3 of Matoussi et al. \cite{MPZ1},  our results can be directly applied to the symmetrical problem, i.e., the RGSDE with an upper barrier. This conclusion is because the proof is only based on a pathwise construction and a fixed-point iteration.
\end{remark}
\subsection{Comparison principle}
\noindent In this subsection, we establish a comparison principle for RGSDEs. At first, we assume additionally a bounded condition on the coefficients $f$, $h$ and $g$ and the obstacle process $S$, and then we remove it in the second step.
\begin{theorem}\label{CT}
Given two RGSDEs that satisfy the assumptions (A1)-(A4), we additionally suppose in the following:
\begin{mylist}
\item $x^1\leq x^2$;
\item $f^i$, $h^i$ and $g^1=g^2=g$ are bounded and $S^i$ are uniformly bounded above, $i=1, 2$;
\item for each $x\in\mathbb{R}$, $f^1_t(x)\leq f^2_t(x)$, $h^1_t(x)\leq h^2_t(x)$; and $S^1_t\leq S^2_t$, $0\leq t\leq T$, q.s..
\end{mylist}
Let $(X^i, K^i)$ be a solution to the RGSDE with data $(f^i, h^i, g, S^i)$, $i=1, 2$, then
\[X^1_t\leq X^2_t,\ 0\leq t\leq T,\ q.s..\]
\end{theorem}
\noindent\textbf{Proof:}
Since $f^i$, $h^i$ and $g$ are bounded and $S^i$ are uniformly bounded above, $i=1, 2$,
using the BDG type inequalities to (\ref{represen}), we deduce that $K^i_T$ has the moment for any arbitrarily large order and for each $t\in[0, T]$, $\lim\limits_{s\rightarrow t}\bar{\mathbb{E}}[|K^i_t-K^i_s|^2]=0$, $i=1, 2$. \\[6pt]
\noindent Notice that $(x^+)^2$ is not a $\mathcal{C}^2(\mathbb{R})$ function. We have to consider $(x^+)^3$ and apply the extended $G$-It\^o's formula to
$((X^1_t-X^2_t)^+)^3$, then
\begin{align}\label{cal1}
((X^1_t-X^2_t)^+)^3
&= 3\int^t_0|(X^1_s-X^2_s)^+|^2(f^1_s(X^1_s)-f^2_s(X^2_s))ds\notag\\
&+ 3\int^t_0|(X^1_s-X^2_s)^+|^2(h^1_s(X^1_s)-h^2_s(X^2_s))d\langle B\rangle_s\notag\\
&+ 3\int^t_0|(X^1_s-X^2_s)^+|^2(g_s(X^1_s)-g_s(X^2_s))dB_s\\
&+ 3\int^t_0|(X^1_s-X^2_s)^+|^2d(K^1_s-K^2_s)\notag\\
&+ 3\int^t_0(X^1_s-X^2_s)^+|g_s(X^1_s)-g_s(X^2_s)|^2 d\langle B\rangle_s.\notag
\end{align}
As on $\{X^1_t>X^2_t\}$, $X^1_t>X^2_t\geq S^2_t\geq S^1_t$, we have
\begin{align}\label{cal2}
\int^t_0|(X^1_s-X^2_s)^+|^2d(K^1_s-K^2_s)&=\int^t_0|(X^1_s-X^2_s)^+|^2dK^1_s-
\int^t_0|(X^1_s-X^2_s)^+|^2dK^2_s\notag\\
&\leq \int^t_0|(X^1_s-S^1_s)^+|^2dK^1_s-
\int^t_0|(X^1_s-X^2_s)^+|^2dK^2_s\\
&\leq -
\int^t_0|(X^1_s-X^2_s)^+|^2dK^2_s\leq 0,\ q.s..\notag
\end{align}
We put (\ref{cal2}) into (\ref{cal1}) and then, by Lipschitz condition (A3) and by taking $G$-expectation on both sides of (\ref{cal2}), we conclude
\begin{align*}\label{complex}
\bar{\mathbb{E}}[((X^1_t-X^2_t)^+)^3]\leq C\bar{\mathbb{E}}\bigg[\int^t_0((X^1_s-X^2_s)^+)^3ds\bigg]\leq C\int^t_0\bar{\mathbb{E}}[((X^1_s-X^2_s)^+)^3]ds.
\end{align*}
Using Gronwall's lemma, it follows that $\bar{\mathbb{E}}[((X^1_t-X^2_t)^+)^3]=0$, which implies the result.\hfill{}$\square$
\begin{theorem}
Given two given RGSDEs that  satisfy the assumptions (A1)-(A4), we additionally suppose in the following:
\begin{mylist}
\item $x^1\leq x^2$ and $g^1=g^2=g$;
\item for each $x\in\mathbb{R}$, $f^1_t(x)\leq f^2_t(x)$ and $h^1_t(x)\leq h^2_t(x)$; and $S^1_t\leq S^2_t$, $0\leq t\leq T$, q.s..
\end{mylist}
Let $(X^i, K^i)$ be a solution to the RGSDE with data $(f^i, h^i, g, S^i)$, $i=1, 2$, then
\[X^1_t\leq X^2_t,\ 0\leq t\leq T,\ q.s..\]
\end{theorem}
\noindent\textbf{Proof:} First, we define the truncated functions for the coefficients and the obstacle process: for each $N>0$, $\xi^N_t(x):=(-N\vee \xi_t(x))\wedge N$, where $\xi$ denote $f^i$, $h^i$, $g$ and $x\in\mathbb{R}$; and $(S^i)^N_t=S^i_t\wedge N$, $0\leq t\leq T$, $i=1, 2$.
It is easy to verify that the truncated coefficients and the obstacle processes satisfy (A2) and (A3).
Moreover, the truncated functions keep the order of the coefficients and the obstacle processes, that is, for each $N>0$,
\[(f^1)^N_t(x)\leq (f^2)^N_t(x)\ {\rm and}\ (h^1)^N_t(x)\leq (h^2)^N_t(x),\ {\rm for\ all}\ x\in\mathbb{R},\ 0\leq t\leq T,\ q.s.;\]
and
\[(S^1)^N_t \leq (S^2)^N_t,\ 0\leq t\leq T,\ q.s..\] 
Consider the following RGSDEs: 
\begin{align*}
(X^i)^N_t=x+\int^t_0 (f^i)^N_s((X^i)^N_s)ds&+\int^t_0 (h^i)^N_s((X^i)^N_s)d\langle B \rangle_s\\ &+\int^t_0 g^N_s((X^i)^N_s) dB_s + (K^i)^N_t,\ 0\leq t\leq T,\ q.s.,\ i=1, 2,
\end{align*}
under the following conditions:\\[6pt]
$
\begin{array}{l}
(a)\ (X^i)^N\in M^p_G([0, T]),\ (X^i)^{N}_t\geq (S^i)^N_t,\ 0\leq t\leq T,\ q.s.;\\[6pt]
(b)\ (K^i)^N\in M_I([0, T])\cap M^p_G([0, T]),\ (K^i)^N_0=0,\ q.s.;\\[6pt]
(c)\ \int^T_0((X^i)^N_t-(S^i)^N_t)d(K^i)^{N}_t=0,\ q.s..
\end{array}
$\\[6pt]
By Theorem \ref{CT}, it is readily observed that for each $N\in\mathbb{N}$, 
\begin{equation}\label{last}
(X^1)^N_t\leq (X^2)^N_t,\ 0\leq t\leq T,\ q.s..
\end{equation}
Meanwhile, by Theorem \ref{345}, we have
\[
\begin{array}{r@{\ }l}
\bar{\mathbb{E}}[\sup\limits_{0\leq s\leq t}|(X^i)^N_s-X^i_s|^p]
&\leq C(\int^T_0(
\bar{\mathbb{E}}[|(f^i)^N_t(X^i_t)-f^i_t(X^i_t)|^p]+\bar{\mathbb{E}}[|(h^i)^N_t(X^i_t)-h^i_t(X^i_t)|^p]\\[6pt]
&+\bar{\mathbb{E}}[|g^N_t(X^i_t)-g_t(X^i_t)|^p]
)dt+\bar{\mathbb{E}}[\sup\limits_{0\leq t\leq T}|(S^i)^N_t-S^i_t|^p]\\[6pt]
&+\int^t_0\bar{\mathbb{E}}[\sup\limits_{0\leq u\leq s}|(X^i)^N_u-X^i_u|^p]ds).
\end{array}
\]
Applying again Gronwall's lemma, we obtain
\begin{align*}
\bar{\mathbb{E}}[\sup_{0\leq t\leq T}|(X^i)^N_t-X^i_t|^p]
&\leq C\bigg(\int^T_0(
\bar{\mathbb{E}}[|(f^i)^N(t, X^i_t)-f^i(t, X^i_t)|^p]+\bar{\mathbb{E}}[|(h^i)^N(t, X^i_t)-h^i(t, X^i_t)|^p]\\
&+\bar{\mathbb{E}}[|g^N(t, X^i_t)-g(t, X^i_t)|^p]
)dt+\bar{\mathbb{E}}[\sup_{0\leq t\leq T}|(S^i)^N_t-S^i_t|^p]\bigg).
\end{align*}
For each $t\in[0, T]$, we calculate
\begin{align*}
\bar{\mathbb{E}}[|(f^i)^N_t(X^i_t)-f^i_t(X^i_t)|^p]&\leq \bar{\mathbb{E}}[|f^i_t(X^i_t)|^p{\mathbb{1}}_{|f^i_t(X^i_t)|>N}]\\
&\leq \bar{\mathbb{E}}[(|f^i_t(0)|+C_L|X^i_t|)^p{\mathbb{1}}_{(|f^i_t(0)|+C_L|X^i_t|)>N}]\\
&\leq C(\bar{\mathbb{E}}[|f^i_t(0)|^p{\mathbb{1}}_{|f^i_t(0)|>\frac{N}{2}}]+\bar{\mathbb{E}}[|X^i_t|^p{\mathbb{1}}_{|X^i_t|>\frac{N}{2}}]).
\end{align*}
Taking into consideration that $f_\cdot(0)$ and $X^i\in M^p_G([0, T])$, from the argument in Remark \ref{remconna}, we know that $f_t(0)$ and $X^i_t\in L^p_G([0, T])$ for almost every $t\in[0, T]$. Therefore, letting $N\rightarrow +\infty$, we have 
\[
\bar{\mathbb{E}}[|(f^i)^N_t(X^i_t)-f^i_t(X^i_t)|^p]\rightarrow 0.
\]
Similarly, we also obtain
\[
\bar{\mathbb{E}}[|(h^i)^N_t(X^i_t)-h^i_t(X^i_t)|^p]\rightarrow 0
;
\]
and
\[
\bar{\mathbb{E}}[|(g^i)^N_t(X^i_t)-g^i_t(X^i_t)|^p]\rightarrow 0.
\]
\noindent Using Lebsgue's dominated convergence theorem to the integrals on $[0, T]$, it follows that
\begin{align}\label{cal3}
\lim_{N\rightarrow +\infty}
\int^T_0(\bar{\mathbb{E}}[|(f^i)^N(t, X^i_t)-f^i(t, X^i_t)|^p]&+\bar{\mathbb{E}}[|(h^i)^N(t, X^i_t)-h^i(t, X^i_t)|^p]\\
&+\bar{\mathbb{E}}[|g^N(t, X^i_t)-g(t, X^i_t)|^p]
)dt=0.\notag
\end{align}
On the other hand, 
\[
\bar{\mathbb{E}}[\sup_{0\leq t\leq T}|(S^i)^N_t-S^i_t|^p]\leq \bar{\mathbb{E}}[\sup_{0\leq t\leq T}(|S^i_t|^p{\mathbb{1}}_{\{|S^i_t|>N\}})]\leq \bar{\mathbb{E}}[\sup_{0\leq t\leq T}|S^i_t|^p{\mathbb{1}}_{\{\sup\limits_{0\leq t\leq T}|S^i_t|>N\}}].
\]
By the proof of Theorem \ref{BXT}, we know that $\sup\limits_{0\leq t\leq T}S^i_t$ is an element in $L^p_G(\Omega_T)$. So we have 
\begin{align}\label{cal4}
\bar{\mathbb{E}}[\sup_{0\leq t\leq T}|(S^i)^N_t-S^i_t|^p]\leq \bar{\mathbb{E}}[\sup_{0\leq t\leq T}|S^i_t|^p{\mathbb{1}}_{\{\sup\limits_{0\leq t\leq T}|S^i_t|>N\}}]\rightarrow 0,\ \textnormal{as}\ N\rightarrow +\infty.
\end{align}
Combining (\ref{cal3}) and (\ref{cal4}), we obtain
\begin{equation}\label{last2}
\bar{\mathbb{E}}[\sup_{0\leq t\leq T}|(X^i)^N_t-X^i_t|^p]\rightarrow 0,\ \textnormal{as}\ N\rightarrow +\infty.
\end{equation}
Then, (\ref{last}) and (\ref{last2}) yield the desired result .\hfill{}$\square$





\ACKNO{The author expresses special thanks to Prof. Ying HU, who provided both the initial inspiration for the work and useful suggestions and also to the anonymous reviewers, who have given me a lot of important constructive advice for the revision.}


\end{document}